# Precise Derivations of Radiative Properties of Porous Media Using Renewal Theory

Shima Hajimirza
Stevens Institute of Technology
Hoboken, NJ, 07030
Email: shajimi1@stevens.edu

## Abstract

This work uses the mathematical machinery of Renewal/Ruin (surplus risk) theory to derive preliminary explicit estimations for the radiative properties of dilute and disperse porous media otherwise only computable accurately with Monte Carlo Ray Tracing (MCRT) simulations. Although random walk and *Lévy* processes have been extensively used for modeling diffuse processes in various transport problems and porous media modeling, relevance to radiation heat transfer is scarce, as opposed to other problems such as probe diffusion and permeability modeling. Furthermore, closed-form derivations that lead to tangible variance reduction in MCRT are widely missing. The particular angle of surplus risk theory provides a richer apparatus to derive directly related quantities. To the best of the authors' knowledge, the current work is the only work relating the surplus risk theory derivations to explicit computations of ray tracing results in porous media. The paper contains mathematical derivations of the radiation heat transfer estimates using the extracted machinery along with proofs and numerical validation using MCRT.

## 1. Introduction

At present, Monte Carlo Ray Tracing (MCRT) simulations are the most reliable non-experimental means for predicting radiative properties of heterogeneous structures[1]–[9]; such predictions are otherwise only computable directly via extremely limited and costly experiments[10]–[16] or approximate models rooted in extensive measurements[17]–[23]. Though only valid in *geometric optics* regimes where particle sizes and interspacing are large (so diffraction and dependent scattering are negligible), MCRT has been persistently invoked and improved over the past few decades [24]–[46]. However, due to stringent precision requirements, the high computational cost of setup and execution, as well as parameter-dependent uncertainty, the required time and computational cost of MC simulations are prohibitively high for obtaining accurate estimates[47]–[50].

The efforts to minimize MCRT cost fall under one of two categories: *statistical characterization* and *variance reduction*. In the face of dispersity and heterogeneity of geometry and material properties, the former attempts to treat the environment as an equivalent *continuous* and *homogenous* scattering domain[51]–[59] and insert the effective radiative statistics into the radiative transfer equations[60]–[62](RTE) solving for all modes of extinction, or as sampling seeds in MCRT. This way the need for full incorporation of the porous structure is obviated[2], [63], [64]. There is a more formal mathematical approach to the statistical treatment of heterogeneous materials called the "homogenization" technique. Rooted in vibration theory, it falls under the general umbrella of *multiscale methods*[65]–[69], though its applicability is limited to highly regular configurations.

Although not always directly intended for MCRT, the theories of mean-beam-length[70]–[72] and the study of configuration or view/shape/angle factors[60], [63], [73]–[75] are within the first category, attempting to make statistical sense of complex geometry by producing equivalent distribution functions. The work along these lines has also studied correction for the co-dependent scattering effects via the use of *dependent factors*[59], [70], [76] which, along with the incorporation of ray transport distance at dispersed media has led to uncommon RTE solutions such as the Dependent Included Discrete Ordinate Method (DI-DOM). For higher accuracy, statistical homogenizing can also be boosted with experimental measurements, *e.g.*, through inverse analysis based on discrete ordinates radiative models[77]. The method of Radiation Distribution Function Identification (RDFI) has been proposed for facilitating MCRT in statistically homogeneous, isotropic porous media modeled as dispersed overlapping spheres[78] and has been extended to real[79]–[81], non-homogeneous and anisotropic beds[82], [83]. In RDFI, the porous medium is characterized by a cumulated extinction coefficient (CEC) distribution function, obtained by minimizing the difference between absorption and extinction distributions rendered by MCRT on a virtual medium and the corresponding experimental quantities. With the advent of computed tomography and imaging techniques, MCRT has been made more accurate and efficient via the so-called upscaling methodology[84], [85], where the characteristic functions of the radiative properties are estimated either through morphological tomography data or from a representative elementary volume (REV) of the medium treated as semi-transparent[20]–[22].

Variance reduction efforts have been more limited and sporadic, and often only applicable when there is prior knowledge of the structure of the radiation solution. The methods are generally either split-based or importance-sampling-based[86]–[88]. Some of the noteworthy methods include the binary spatial partitioning (BSP) grouping algorithm[89], probabilistic multiple-rays tracing technique[90], and the delta-scattering method[88], [91]. Analytical methods are also effective in reducing the variance of MCRT by providing rapid initial estimates. Albeit limited and nonuniversal, example work along those lines are those based on the Spatial Averaging Theorem (SAT)[92]–[94], as well as custom analytical and empirical models constructed from MCRT ground truth data[24], [25], [76], [95]–[97], though the scopes are limited to specific packing distributions and/or material compositions. Recent data-driven endeavors to use machine learning for MCRT estimations[98]–[104] can also be viewed as unorthodox branches of variance reduction efforts, though the work is in its very infancy.

The present work offers a novel approach to variance reduction based on the mathematics of Renewal/Ruin theory (Cramér-Lundberg), providing more generalizable and accurate results than the existing variance reduction schemes. The essence of the models is in approximating the travelling ray within a porous bed with a one- or two-sided renewal random walk process. Using the existing machinery of surplus risk theory and ruin analysis, we can then derive close upper bounds and explicit tight estimates for radiative metrics. Although random walk and *Lévy* processes have been extensively used for modeling diffuse processes[105]–[109] in various transport problems and porous media modeling, there are two shortcomings: (i) relevance to radiation heat transfer is scarce, as opposed to other problems such as probe diffusion and permeability modeling[110]–[113], (ii) closed-form derivations that lead to tangible variance reduction in MCRT are missing. The particular angle of surplus risk theory provides a richer apparatus to derive directly related quantities. To the best of the authors knowledge, the current work is the only work relating the surplus risk theory derivations to explicit computations of ray tracing results in porous media.

**Summary of Results**

We leverage key theorems of Ruin theory based on Markovian processes[114] to derive explicit characterizations of power extinction and reflectivity as would be derived by MCRT in flat one- and two-sided polydisperse porous media (Fig. 1). Considering a traveling beam that enters a flat (Beerian) porous medium and experiences independent mean free path and subsequent reflection events with known statistics $F_\ell(\cdot)$ (cumulative distribution function), the objective is to derive reflection, absorption, and transmission fractions. Assuming the media is an infinitely horizontally stretched slab, the sequence of vertical position (depth) of the ray at the time of scattering events can be modeled as a 1D Markov chain random walk. Existing results from Ruin theory (*e.g.*, [114]) provide explicit calculations for the joint moment generating function of Ruin time random variables, including stopping time (see Lemmas 2 and 4 in Section 3 for a summary of relevant extractions). By invoking a series of techniques including Jensen[115] and Cauchy-Schwarz inequalities[116], Delta method[117], and Wald's identity[118], we successfully relate the moment generating functions to approximate collective power dissipation and reflection in a statistical Beerian and non-Beerian porous media, both for one-sided and two-sided beds (See Theorems 2 and 3). We compare the derived analytical results with those obtained from exhaustive Monte Carlo ray tracing simulations and report astonishingly close matches, serving as validation for our theoretical results.

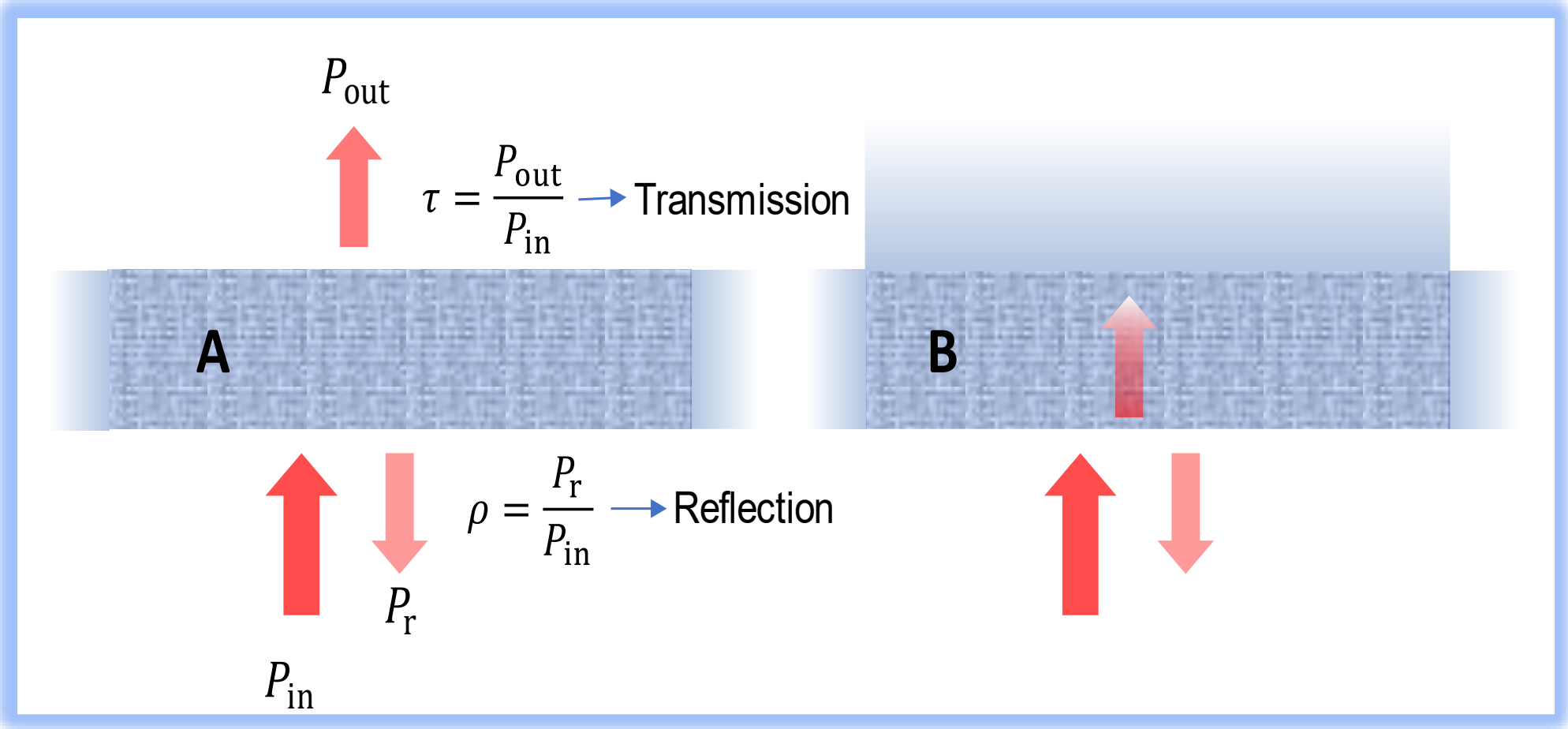

**Figure 1.** (A) two-sided (thick) and (B) one-sided infinite porous slab.

## 2. Mathematical Derivations

### 2.1. Part I. One-sided medium.

**Definition 1**. A one-sided $(x, p, F_y)$-renewal process (See Figure 2) is a random walk that starts at point $x \geq 0$, and takes consecutive steps of sizes $Y_i, i > 0$ with cumulative distribution $F_y$, where the step is to the right (increasing) with probability $p$ and to the left (decreasing) with probability $q = 1 - p$, and terminates as soon as it becomes less than 0. $T_x$ and $Z_x$ are defined as random variables describing the stopping (exit) time and overshoot at the stopping time respectively. In other words:

$$T_x = \min\{T \geq 0 \mid \sum_{i=1}^{T} s_i Y_i < 0\} \tag{1}$$

$$Z_x = -\sum_{i=1}^{T_x} s_i Y_i - x, \tag{2}$$

where $s_i \in \{1, -1\}$ is the side of the $i^{\text{th}}$ step.

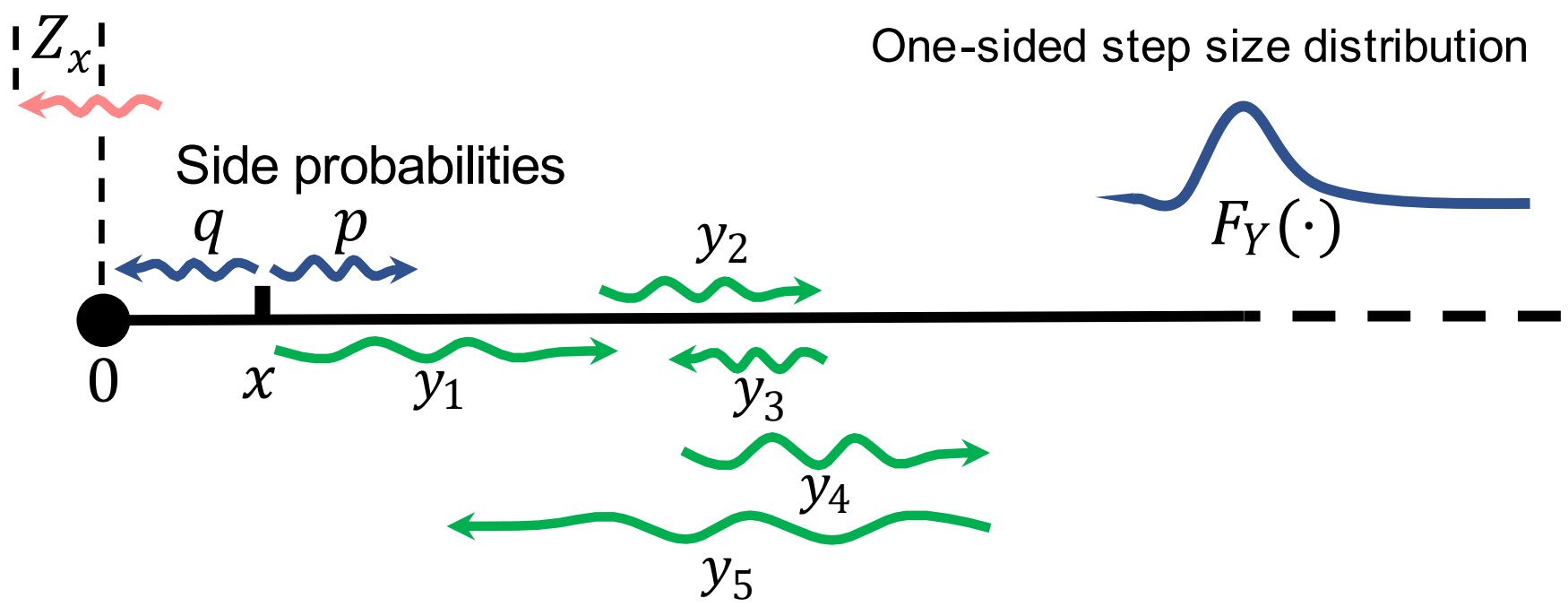

**Figure 2.** One-sided renewal process.

**Lemma 1.** In a one-sided flat polydisperse porous medium with dissipation factor $\beta$ and opaque particles, the following upper bound always holds for the reflected (non-absorbed) fraction $\rho$ of the incident power:

$$\rho \leq \rho^u = \mathbb{E}_x\left(e^{-\beta x}\sqrt{\mathbb{E}_{T_x}\left(\mathbb{E}_Y\left(e^{-2\beta Y}\right)\right)^{T_x}\mathbb{E}_{Z_x}e^{2\beta|Z_x|}}\right), \tag{3}$$

where $x$ is the depth of the first scattering location inside the medium, and $T_x, Z_x$ are the stopping time and overshoot random variables associated with the one-sided $(x, 0.5, F_y)$-renewal process constructed from the depth series of the consecutive scattering (reflecting) locations (see Figure 3).

Proof. Assuming that scattering is isotropic, there is an equal probability of decreases versus increases in the beam depth, *i.e.*, $p = q = 0.5$. Therefore, considering that the medium has infinite height, the sequence of ray depth values at consecutive scattering locations is a one-sided $(x, 0.5, F_y)$ renewal process. Denoting by $\rho_x$ the reflected power of a beam of unity power that undergoes an isotropic scattering at a depth $x$ inside the bed, the average reflect power is the expected value of $e^{-\beta x}\rho_x$ for all possible initial scattering depth values $x$. On the other hand, $\rho_x$ can be expressed in terms of renewal process variables. Firstly, it can be calculated from the Beer law expressed as $\rho_x = e^{-\beta L_x + \beta|Z_x|}$, where $L_x$ and $Z_x$ are the total random walk travel and exit overshoot by the stopping time, respectively. Furthermore, $L_x = \sum_i^{T_x} y_i$ where $T$ is the stopping time and $y_i$s are absolute changes in the depth of the ray between every two consecutive scattering events. We can now invoke Cauchy-Schwarz inequality to conclude that:

$$\mathbb{E}_x\rho_x = \mathbb{E}_x\left(e^{-\beta x}\mathbb{E}\left(e^{-\beta L_x+\beta|Z_x|}\right)\right) \le \mathbb{E}_x\left(e^{-\beta x}\sqrt{\mathbb{E}e^{-2\beta L_x}\mathbb{E}e^{2\beta|Z_x|}}\right).$$

Furthermore, from Wald's identity we have $\mathbb{E}e^{-2\beta L_x} = \mathbb{E}e^{-2\beta\sum_i^{T_x} y_i} = \left(\mathbb{E}_Y\left(e^{-2\beta Y}\right)\right)^{T_x}$ finalizing the proof□

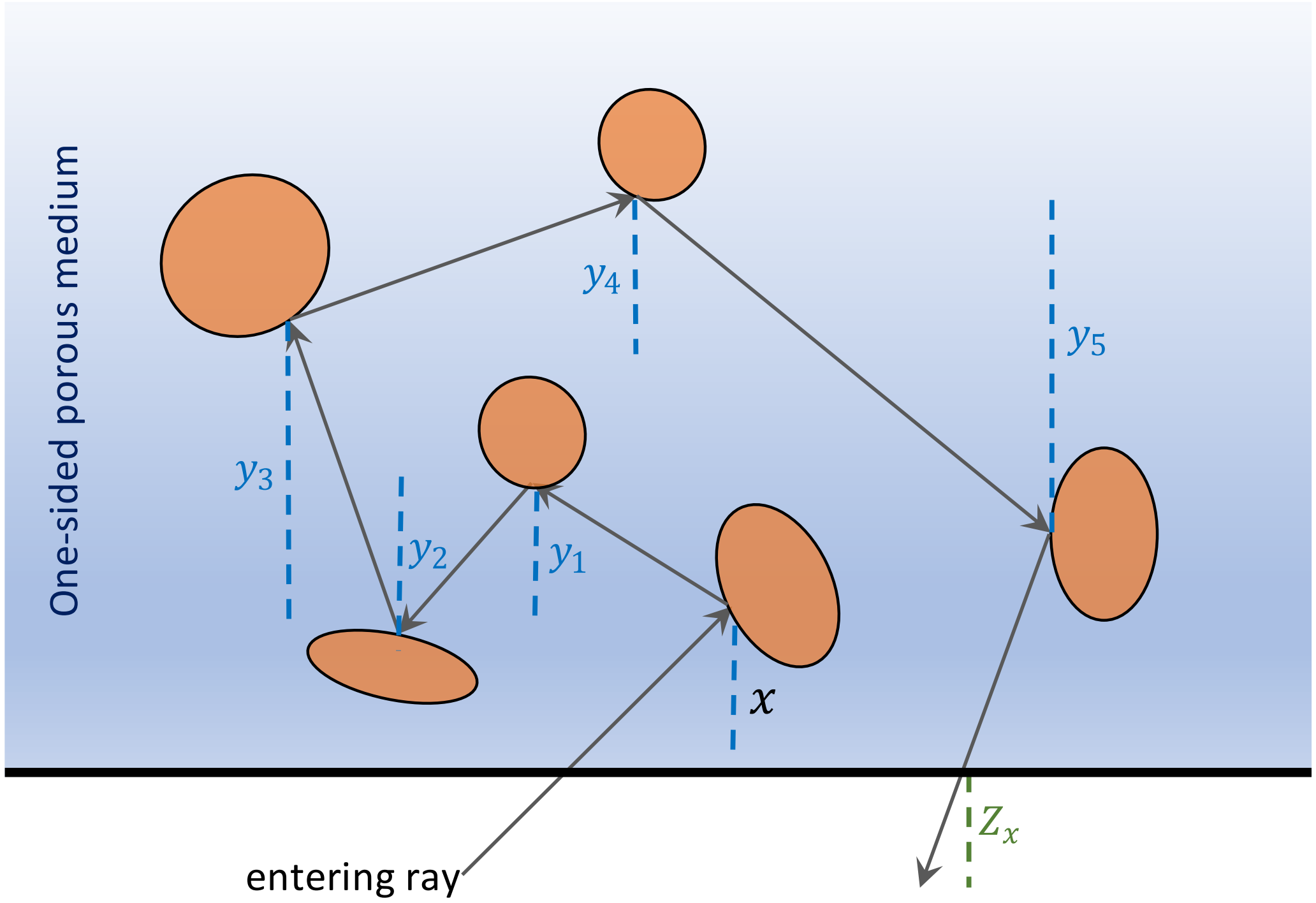

**Figure 3.** Schematics of ray travel in a one-sided porous media, with consecutive scattering events. The ray travel can be modeled with a one-sided renewal process representing the depth (vertical) location of the ray at locations where it hits opaque particles.

**Theorem 1**. In a one-sided flat polydisperse porous medium with dissipation factor $\beta$ and opaque particles, the following approximation holds for the reflected (non-absorbed) fraction $\rho$ of the incident power:

$$\rho \approx \hat{\rho} = \left(e^{-\beta x}\mathbb{E}_{T_x}\mathbb{E}_Y\left(e^{-2\beta Y}\right)^{T_x}\right)\left(1 - \Delta_\epsilon \mathbb{E}_{Z_x} e^{-\epsilon Z_x/\bar{x}}\,\beta\right). \tag{4}$$

Proof. Directly follows from the use of Delta method:

$$Z_m = \Delta_\epsilon \mathbb{E}e^{-\epsilon/\bar{x}\mu Z}$$

$$\mathbb{E}e^{\beta Z_x} \approx 1 - \Delta_\epsilon \mathbb{E}e^{-\epsilon Z/\bar{x}}\,\beta$$

$$\rho = \mathbb{E}_x\left(e^{-\beta x}\mathbb{E}\left(e^{-\beta L_x+\beta|Z_x|}\right)\right) \approx \mathbb{E}\left(e^{-\beta x}\mathbb{E}e^{-\beta L_x}\right)\left(1-\Delta_\epsilon \mathbb{E}e^{-\frac{\epsilon Z}{\bar{x}}}\beta\right) \approx \hat{\rho}$$
$$= \mathbb{E}\left(e^{-\beta x}\mathbb{E}\left(\mathbb{E}e^{-2\beta Y}\right)^{T_x}\right)\left(1-\Delta_\epsilon \mathbb{E}e^{-\frac{\epsilon Z}{\bar{x}}}\beta\right).$$

□

**Lemma 2** (Combination of results from [114] for one-sided renewal process). The following explicit formula holds for the moment generating function of the overshoot $Z_x$ and stopping time $T_x$ of a one-sided $(x, p, F_y)$-renewal process:

$$\mathbb{E}_{T_x Z_x}\left(\alpha^{T_x} e^{-\zeta Z_x}\right) = \sum_{i=1}^{m} c_i e^{\gamma_i x}, \forall 0 \leq \alpha \leq 1, \zeta \geq 0, \tag{5}$$

where $\gamma_i$s are the non-positive roots of the Cramer-Lundberg equation:

$$pL_+(\gamma) + qL_-(\gamma) = \alpha^{-1},$$

where $L_+$ and $L_-$ are one-sided Laplace transforms of the CDF $F_y$:

$$L_+(v \in \mathbb{C}^+) = \mathcal{L}_{F_y}^+(v) = \frac{P_+(v)}{R_+(v)},$$

$$L_-(v \in \mathbb{C}^+) = \mathcal{L}_{F_y}^-(v) = \frac{P_-(v)}{R_-(v)},$$

where $P_+, P_-, R_+, R_-$ are polynomials of degree $m$. The coefficients $c_i$s are given by:

$$c_i = \frac{R_+(\gamma_i)}{R_+(\zeta)} \frac{\Pi_{j \neq i}(\zeta - \gamma_j)}{\Pi_{j \neq i}(\gamma_i - \gamma_j)}.$$

Figure 4 summarized these relations in a compact way for convenience.

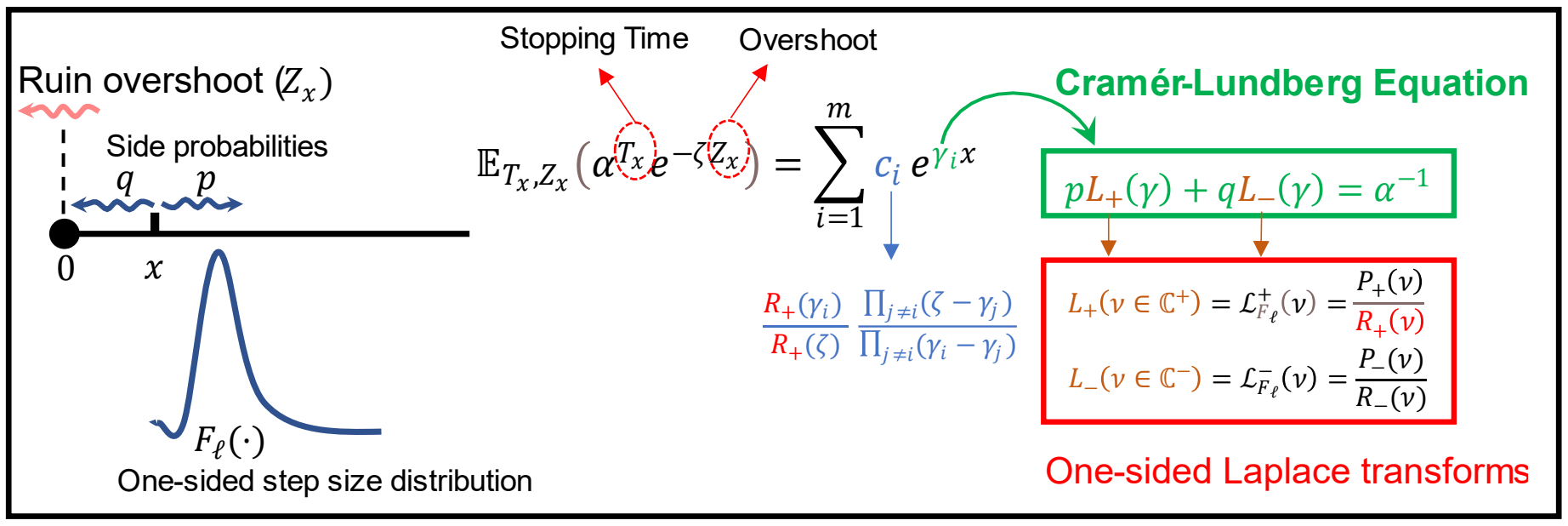

**Figure 4.** Summary of mathematical relations of Lemma 2 for a one-sided renewal process.

**Theorem 2.** If the homogenized beam length distribution in a Beerian porous media with opaque particles and dissipation factor $\beta$ can be expressed as a single exponential distribution with the rate parameter $\mu$, *i.e.*, $F_\ell(y) = \mu e^{-\mu y}$, then the upper bound ($\rho^u$) and approximation $\hat{\rho}$ for the normal reflection become:

$$\hat{\rho} = \frac{1-\sqrt{\frac{2\eta}{2\eta+1}}}{\eta+1+\sqrt{\frac{2\eta}{2\eta+1}}}(1+\eta) \approx \rho \leq \rho^u = \frac{\left(1-\sqrt{\frac{2\eta}{1+2\eta}}\right)^{\frac{1}{2}} \left(\frac{1}{1-2\eta}\right)^{\frac{1}{2}}}{1+\eta+0.5\sqrt{\frac{2\eta}{1+2\eta}}}, \tag{6}$$

where $\eta = \beta/\mu$. The generalization to non-normal incident angle $\theta$ is as follows:

$$\hat{\rho} = \left(1 - \sqrt{\tfrac{2\eta}{1+2\eta}}\right)\left(1 + \left(\eta + \sqrt{\tfrac{2\eta}{1+2\eta}}\right)\cos\theta\right)^{-1}(1+\eta). \tag{7}$$

Proof: Given the single exponential distribution assumption, (5) simplifies to $\mathbb{E}_{T_x,Z_x}\left(\alpha^{T_x}e^{-\zeta Z_x}\right) = ce^{\gamma x}$ where for the case of single mode distribution:

$$c = \frac{R_+(\gamma)}{R_+(\zeta)} = \frac{\gamma+\mu}{\zeta+\mu}$$

$$\gamma = -\mu\sqrt{1-\alpha}, \qquad \alpha = \frac{\mu}{\mu+2\beta}$$

First, we prove the left-hand side of (6). We have from (4) that:

$$\hat{\rho} = \mathbb{E}\left(e^{-\beta x}\mathbb{E}\alpha^{T_x}\right)\left(1 - \Delta_\epsilon \mathbb{E}e^{-\epsilon Z/\bar{x}}\,\beta\right)$$

where $\alpha = \mathbb{E}_Y\left(e^{-2\beta Y}\right)$, $\bar{x} = \mathbb{E}x$. For the case of exponential distribution with parameter $\mu$, $\alpha = \frac{\gamma}{\beta+\gamma}$ and $\bar{x} = 1/\mu$, and therefore we can write:

$$-\Delta_\epsilon \mathbb{E}e^{-(\epsilon/\bar{x})Z} = -\frac{c_1\left(\frac{2\epsilon}{\bar{x}}\right) - c_1\left(\frac{\epsilon}{\bar{x}}\right)}{\frac{\epsilon}{\bar{x}}} = -\frac{\frac{\mu}{\frac{2\epsilon}{\bar{x}}+\mu} - \frac{\mu}{\frac{\epsilon}{\bar{x}}+\mu}}{\frac{\epsilon}{\bar{x}}} = \frac{\frac{1}{\mu}}{(2\epsilon+1)(\epsilon+1)} \approx \frac{1}{\mu}.$$

Furthermore, from the moment generating results of Lemma 2 we know that:

$$\mathbb{E}\alpha^{T_x} = c_1 e^{\gamma x},$$

where

$$c_1(\zeta = 0) = \frac{R_+(\gamma)}{R_+(\zeta=0)} = \frac{\gamma+\mu}{\mu},$$

$$\gamma = -\mu\sqrt{1-\alpha} = -\mu\sqrt{\frac{2\beta}{2\beta+\mu}}.$$

Therefore:

$$\hat{\rho} = \mathbb{E}\left(c_1(0)e^{-(\beta-\gamma)x}\right)\left(1+\frac{\beta}{\mu}\right) = \frac{\gamma+\mu}{\beta+\mu-\gamma}\left(1+\frac{\beta}{\mu}\right) = \frac{-\mu\sqrt{\frac{2\beta}{2\beta+\mu}}+\mu}{\beta+\mu+\mu\sqrt{\frac{2\beta}{2\beta+\mu}}}\left(1+\frac{\beta}{\mu}\right)$$

$$= \frac{1-\sqrt{\frac{2\eta}{2\eta+1}}}{\eta+1+\sqrt{\frac{2\eta}{2\eta+1}}}(1+\eta).$$

The case of non-zero incident angle is a very straightforward modification by considering:

$$\rho_\theta = \mathbb{E}_x \rho_{\tilde{x}},$$

where $\tilde{x} = x\cos\theta$. Therefore:

$$\rho_\theta = \int_0^\infty \mu e^{-(\mu+\beta\cos\theta)x}\ \mathbb{E}\left(e^{-\beta L_{\tilde{x}}+\beta|Z_{\tilde{x}}|}\right)dx.$$

The tight estimate $\hat{\rho}_\theta$ can therefore be written as:

$$\hat{\rho}_\theta = \mathbb{E}\left(e^{-\beta\tilde{x}}\mathbb{E}\alpha^{T_{\tilde{x}}}\right)\left(1 - \Delta_\epsilon \mathbb{E}e^{-\epsilon Z/\mathbb{E}\tilde{x}}\beta\right).$$

where $\alpha = \mathbb{E}_Y\left(e^{-\beta Y}\right)$. Using similar calculations as before, we can write:

$$\hat{\rho} = \mathbb{E}\left(c_1(0)e^{-(\beta-\gamma)x\cos\theta}\right)\left(1+\frac{\beta}{\mu}\right) = \frac{\gamma+\mu}{\mu+(\beta-\gamma)\cos\theta}\left(1+\frac{\beta}{\mu}\right)$$

$$= \frac{1-\sqrt{\frac{2\eta}{1+2\eta}}}{1+\left(\eta+\sqrt{\frac{2\eta}{1+2\eta}}\right)\cos\theta}(1+\eta),$$

which proves (5). We now turn to the right-hand side of (4). We can write:

$$\rho^{(u)} = \mathbb{E}_x\left(e^{-\beta x}\sqrt{\mathbb{E}\alpha^{T_x}\mathbb{E}_{Z_x}e^{2\beta|Z_x|}}\right).$$

Again, as shown above, $\mathbb{E}\alpha^{T_x} = \frac{\gamma+\mu}{\mu}e^{\gamma x}$. Furthermore, we can use the moment-generating function of $Z_x$ from Lemma 2 to conclude that $Z_x$ has an exponential distribution with rate parameter $\mu$. Therefore, (note that it is only valid for $\mu > 2\beta$):

$$\mathbb{E}e^{2\beta Z_x} = \frac{\mu}{\mu-2\beta}.$$

Putting all these together along with $\gamma = -\mu\sqrt{\frac{2\beta}{2\beta+\mu}}$ we get:

$$\rho^{(u)} = \mathbb{E}_x\left(e^{-\beta x}\sqrt{e^{\gamma x}\left(\frac{\gamma+\mu}{\mu}\right)\left(\frac{\mu}{\mu-2\beta}\right)}\right) = \kappa\mathbb{E}e^{-x(\beta-\gamma/2)} = \frac{\kappa\mu}{\mu+\beta+\gamma/2},$$

where $\kappa = \sqrt{\left(\frac{\gamma+\mu}{\mu}\right)\left(\frac{\mu}{\mu-2\beta}\right)}$, which can be simplified as:

$$\rho^{(u)} = \sqrt{\left(\frac{\gamma+\mu}{\mu}\right)\left(\frac{\mu}{\mu-2\beta}\right)}\frac{\mu}{\mu+\beta-\gamma}$$

$$= \left(1 - \sqrt{\frac{2\beta}{2\beta + \mu}}\right)^{0.5} \left(\frac{\mu}{\mu - 2\beta}\right)^{0.5} \frac{\mu}{\mu + \beta + 0.5\mu\sqrt{\frac{2\beta}{2\beta + \mu}}}$$

$$= \left(1 - \sqrt{\frac{2\eta}{2\eta + 1}}\right)^{0.5} \left(\frac{1}{1 - 2\eta}\right)^{0.2} \frac{1}{1 + \eta + 0.5\sqrt{\frac{2\eta}{2\eta + 1}}},$$

thus completing the proof □

**Remark.** Note that the dissipation factor in a Beerian environment is given by $\beta = \frac{2\pi k}{\lambda}$, where $k$ is the complex component of the refractive index. Therefore, the upper bound (right hand side) in Theorem 2 is only valid for

$$\mu > 2\beta \rightarrow \mu > \frac{4\pi k}{\lambda} \rightarrow \lambda > 4\pi kL,$$

where $L$ is the average mean free beam length penetration length $\mathbb{E}\ell, \ell \sim \text{Exp}(\mu)$. So this means the upper bound is only valid in microscale porous media in which the average distance to closest scattering particle within the participating medium is comparable to wavelength (near field).

**2.2. Part II. Two-sided medium.**

**Definition 2.** A two-sided $(x, h, p, F_y)$-renewal process (See Fig. 5) is a random walk that starts at point $x \geq 0$, and takes consecutive steps of sizes $Y_i, i \geq 0$ with cumulative distribution $F_y$, where steps are to the right (increasing) with probability $p$ and to the left with probability $q = 1 - p$, and terminates as soon as it becomes less than 0 or larger than $h$. $T_{x,h}$, $Z^+_{x,h}$, $Z^-_{x,h}$ are defined as random variables describing the stopping time and positive and negative overshoots at stopping time, respectively. In other words:

$$T_x = \min\{T \geq 0 \mid \sum_{i=1}^{T} s_i Y_i < 0 \text{ or } \sum_{i=1}^{T} s_i Y_i > h\}, \tag{7}$$

$$Z^+_{x,h} = \max\left(\textstyle\sum_{i=1}^{T_x} s_i Y_i - h, 0\right), \tag{8}$$

$$Z^-_{x,h} = \max\left(-\sum_{i=1}^{T_x} s_i Y_i - x, 0\right), \tag{9}$$

where $s_i \in \{1, -1\}$ is the side of the $i^{\text{th}}$ step.

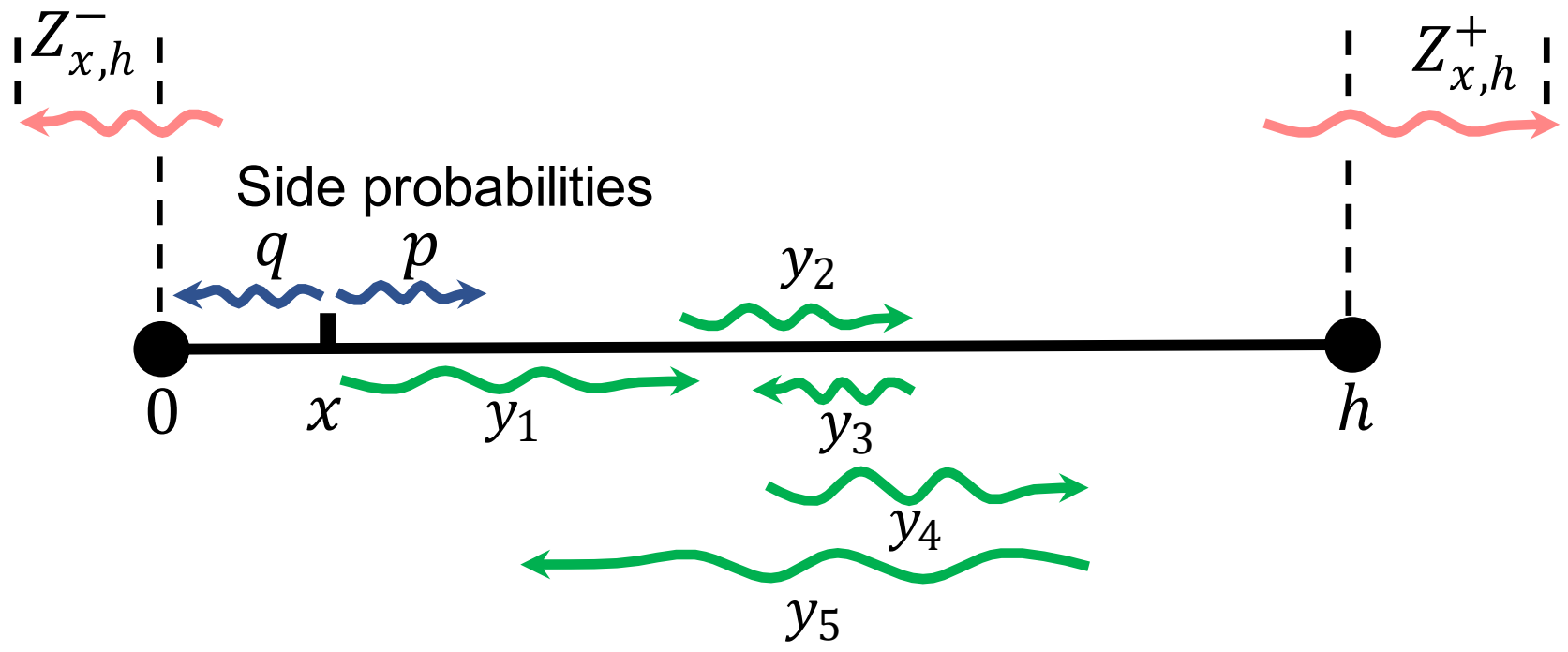

**Figure 5.** Two-sided renewal process.

**Lemma 3.** In a two-sided flat polydisperse non-Beerian porous bed of height $h$, the reflectivity fraction $\rho_\theta$ of the incident power at angle $\theta$ is given by:

$$\mathbb{E}_x \rho_{\tilde{x}} = \mathbb{P}(X_\infty \leq 0) = \frac{h - x + \mathbb{E}Z^+_{h,\tilde{x}}}{\mathbb{E}Z^-_{h,\tilde{x}} + \mathbb{E}Z^+_{h,\tilde{x}} + h}, \tag{10}$$

where $x$ is a random variable representing the first scattering depth inside the medium and $\tilde{x} = x \cos \theta$, and $Z^+_{h,\tilde{x}}, Z^-_{h,\tilde{x}}$ are the overshoot variables defined for the equivalent two-sided renewal process of Definition 2.

**Lemma 4** (Combination of results from [114] for two-sided renewal processes). The following explicit formula holds for the moment generating function of the overshoots $Z^+_{x,h}, Z^-_{x,h}$ and stopping time $T_{x,h}$ of a two-sided $(x, h, p, F_y)$-renewal process:

$$\mathbb{E}\left(A\alpha^{T_{x,h}} e^{-\zeta Z^-_{x,h}}\right) + \mathbb{E}\left(B\alpha^{T_{x,h}} e^{-\xi Z^+_{x,h}}\right) = \sum_{i=1}^{m} c_i e^{\gamma_i x} + \sum_{i=1}^{m} d_i e^{\delta_i x}, \quad \forall 0 \leq \alpha \leq 1, \zeta, \xi \geq 0, A, B \in \mathbb{R} \tag{11}$$

where $\gamma_i$s and $\delta_i$s are the non-positive roots of the Cramer-Lundberg equation:

$$pL_+(\gamma) + qL_-(\gamma) = \alpha^{-1},$$

where $L_+$ and $L_-$ are one-sided Laplace transforms of the CDF $F_y$:

$$L_+(v \in \mathbb{C}^+) = \mathcal{L}^+_{F_y}(v) = \sum_{i=1}^{m} a_i \frac{\mu_i}{\mu_i + v},$$

$$L_-(v \in \mathbb{C}^+) = \mathcal{L}^-_{F_y}(v) = \sum_{i=1}^{l} b_i \frac{\psi_i}{\psi_i + v}.$$

The coefficients $c_i$s and $d_i$s are the solutions of the following linear system of equations:

$$\sum_{i=1}^{m} \frac{c_i}{\mu_t + \gamma_i} + \sum_{i=1}^{l} \frac{d_i}{\mu_t + \delta_i} = \frac{A}{\mu_t + \zeta}, \qquad 1 \leq t \leq m$$

$$\sum_{i=1}^{m} \frac{c_i e^{\gamma_i h}}{\psi_t - \gamma_i} + \sum_{i=1}^{l} \frac{d_i e^{\delta_i h}}{\psi_t - \delta_i} = \frac{B}{\psi_t + \xi}, \qquad 1 \le t \le l$$

Fig. 6 summarized these relations in a compact way for convenience.

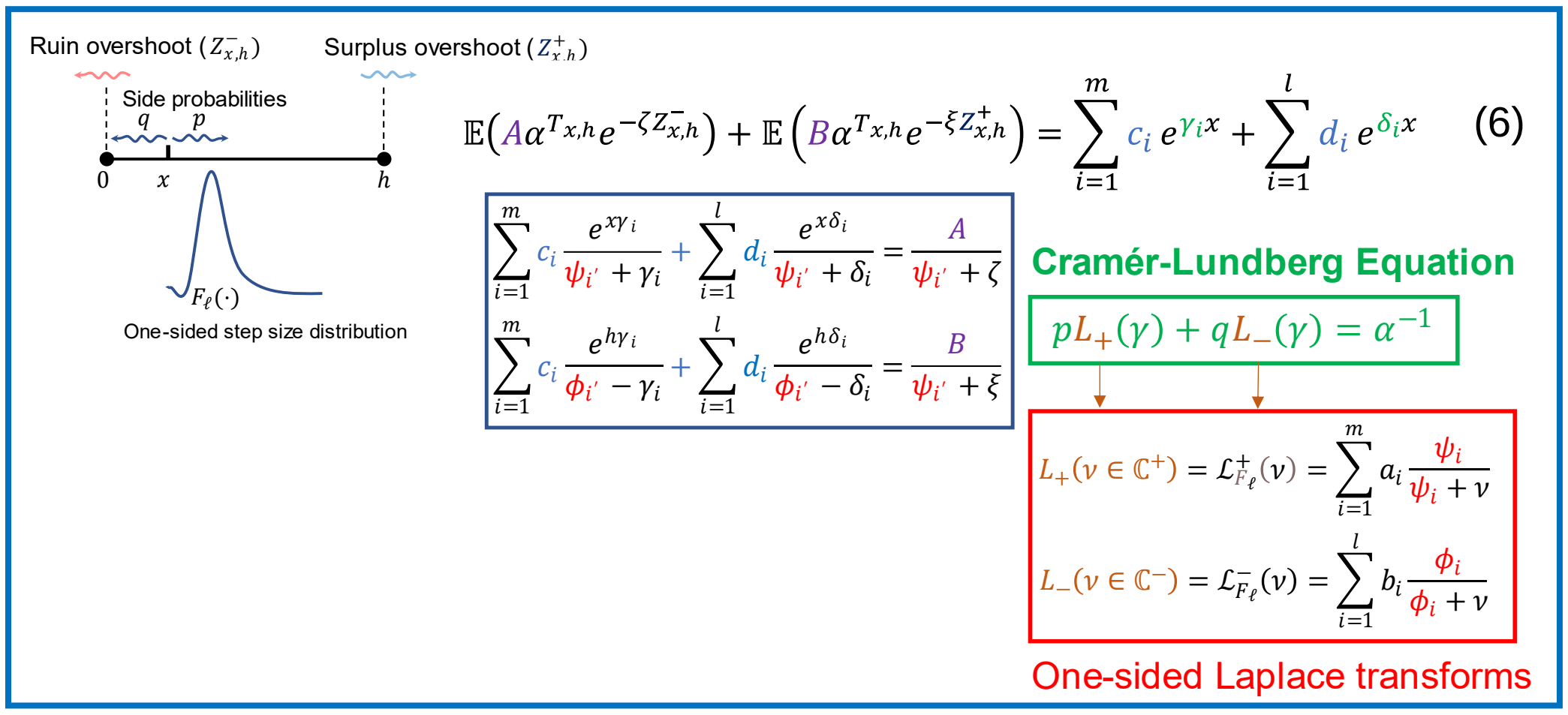

**Figure 6.** Summary of mathematical relations of Lemma 2 for a one-sided renewal process.

**Theorem 3.** If the homogenized beam length distribution in a non-Beerian two-sided porous bed of height $h$ with opaque particles can be expressed as a single exponential distribution with rate parameter $\mu$, i.e., $F_\ell(y) = \mu e^{-\mu y}$, then the reflection fraction $\rho_\theta$ for incident angle $\theta$ is <u>precisely</u> equal to:

$$\rho_\theta = \frac{(1-\cos\theta)\left(1-e^{-\frac{h\mu}{\cos\theta}}\right)+h\mu}{h\mu+2}. \tag{12}$$

Proof. We can easily follow Lemma 4 to conclude that $\mathbb{E}Z^+_{x,h} = Z^-_{x,h} = 1/\mu$. Using that directly in Lemma 3 results in:

$$\mathbb{E}\rho_{x\cos\theta} = \mathbb{E}_x\left(\frac{h - x\cos\theta + 1/\mu}{2/\mu + h}\middle| x\cos\theta \le h\right).$$

The rest of the proof follows mechanically from calculating the conditional expected value given the distribution of $x$□

## 3. Numerical Results

The analytical upper bound and estimate values of Eqn. (4) are compared with those derived from exhaustive Monte Carlo simulations in Fig. 7, assuming a statistical Beerian ($\beta = 1$) infinite medium with an exponential mean free beam length distribution with parameter $\mu$. Two regimes of $\mu$ are considered, namely near- and far-field (Fig.s 7.A and 7.B). The reflectivity estimate $\hat{\rho}$ is accurate down to an astonishing 0.01%, as it is directly stemming from analytical derivations.

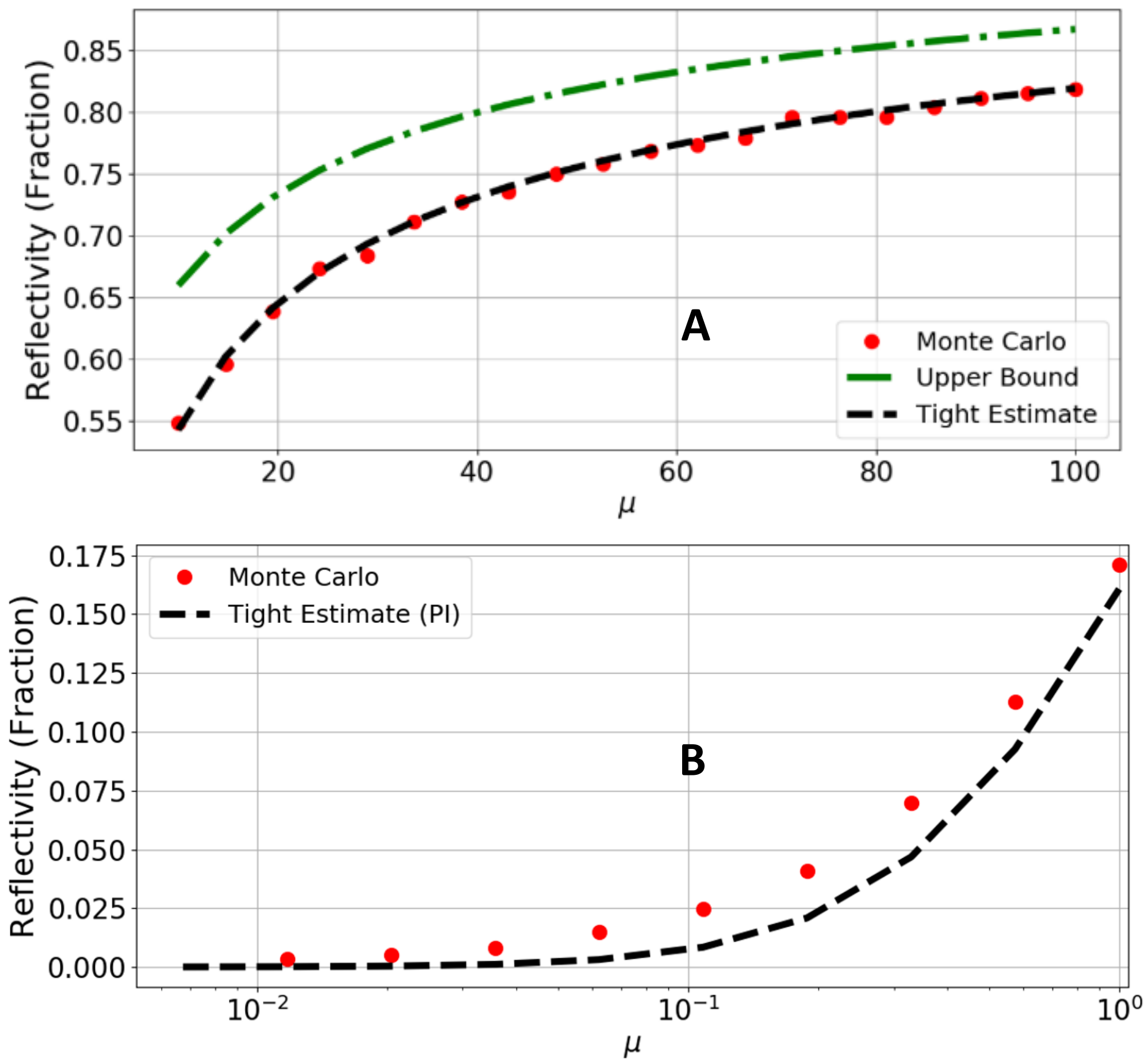

**Figure 7**. Comparison of the analytical reflectivity upper bound and delta-method estimate of the PI with exhaustive MC simulations for a Beerian one-sided porous medium. (A) near field, (B) far-field.

To demonstrates the efficacy of using statistical homogenization in deriving approximate radiative properties, we studied a one-sided Beerian packed bed (infinite slab, $\beta > 0$) filled with random overlapping circular particles (Figure 8), assuming 100% opacity, and computed the empirical distribution of mean free beam length using an in-house pixel-based MCRT code. Fig. 9 shows the average radiation flux for the normal incident angle inside the medium obtained using MCRT using 10000 iterations. The empirical distribution of the mean free beam length, $f_\ell(x) = \partial F_\ell / \partial x$, namely, the random variable describing the distance between every two consecutive scattering events was then fit using a single-rate exponential distribution with Maximum-Likelihood and least-square PDF fitting. The empirical and analytical distributions are shown in Fig. 10. An analytical estimate of the reflection of upward incident radiation was obtained by introducing the fitted distribution $\hat{F}_\ell(\cdot)$ to the results of Theorem 2 and is compared to computational values from MCRT and pack-free MC using 10000 iterations, as depicted in Fig. 11. The curves indicate a high estimation accuracy for the purely analytical model.

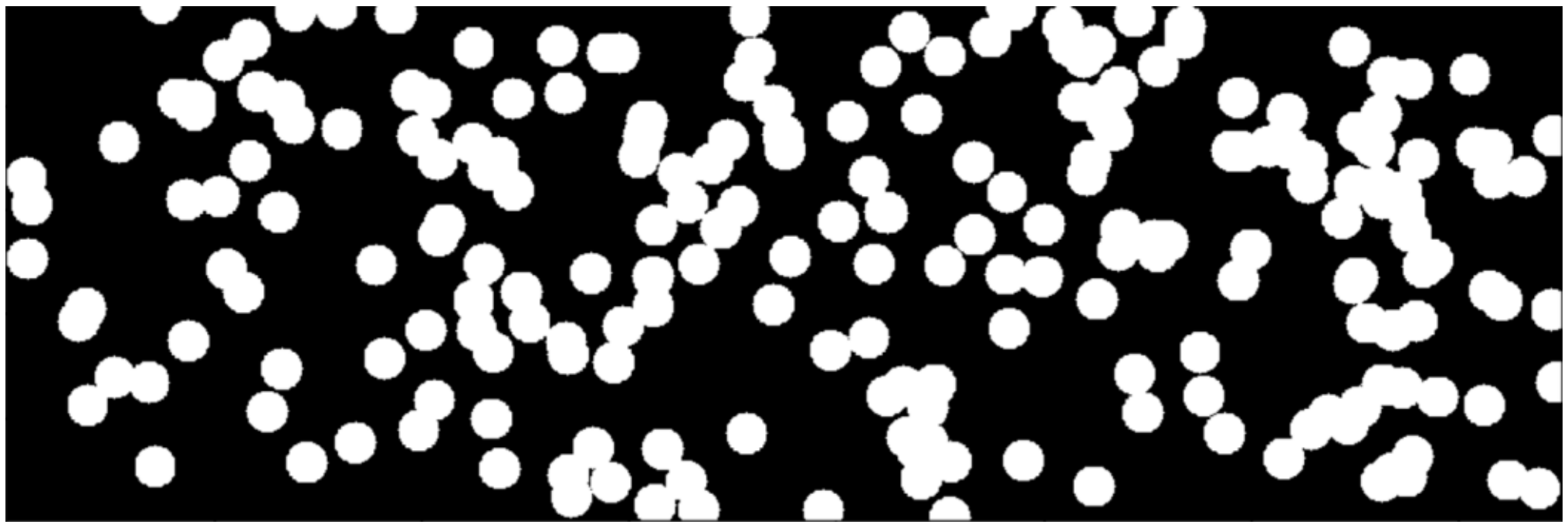

**Figure 8.** Sample random one-sided porous medium filled with overlapping circular opaque particles.

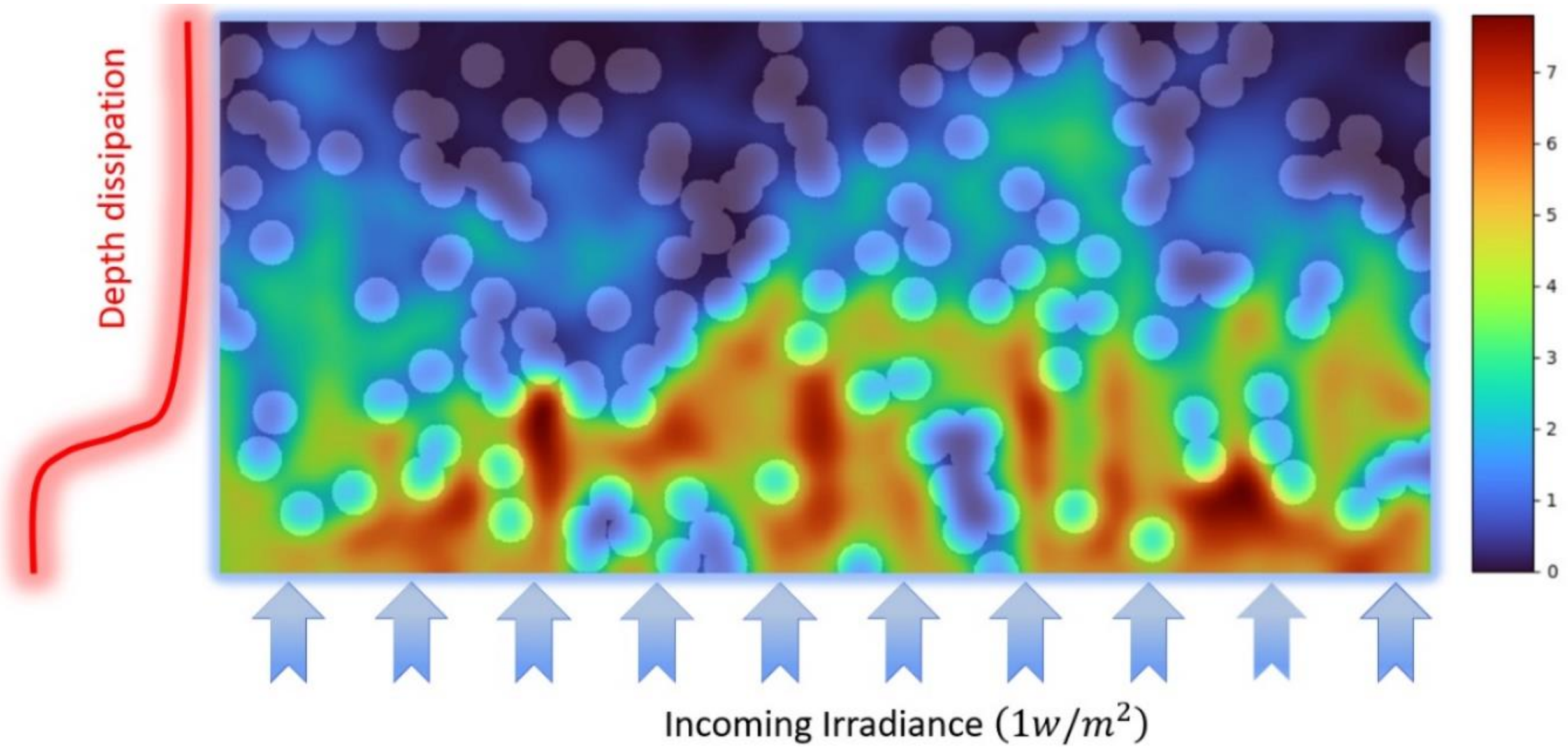

**Figure 9.** Average radiation flux calculated via ray tracing simulations with 10000 iterations for the normal incident (from the bottom side) for the one-sided medium of Figure 8.

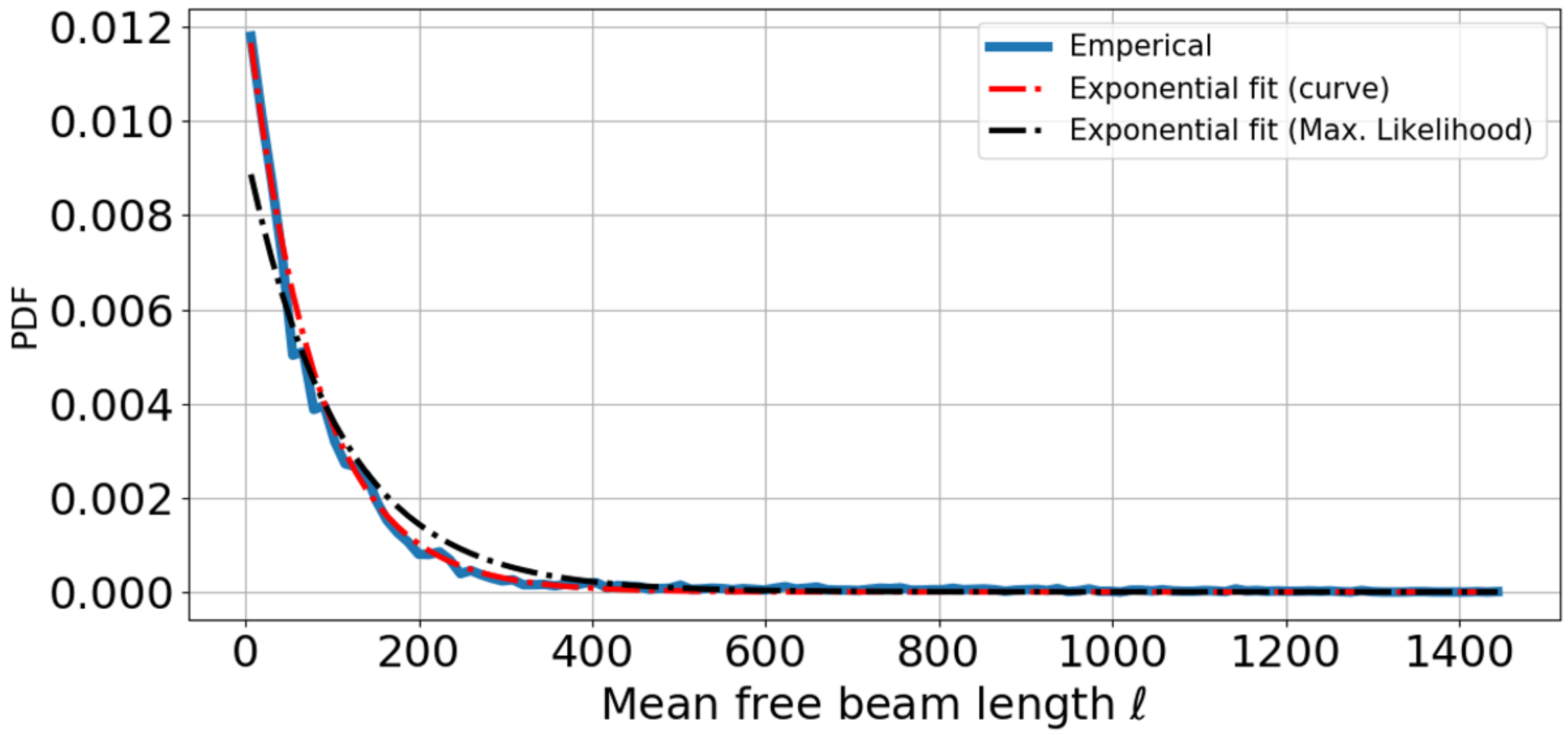

**Figure 10.** Empirical homogenized probability distribution function (PDF) of the mean free beam length inside the porous medium of Figure 8, along with exponential approximations obtained via curve fitting and maximum likelihood approaches.

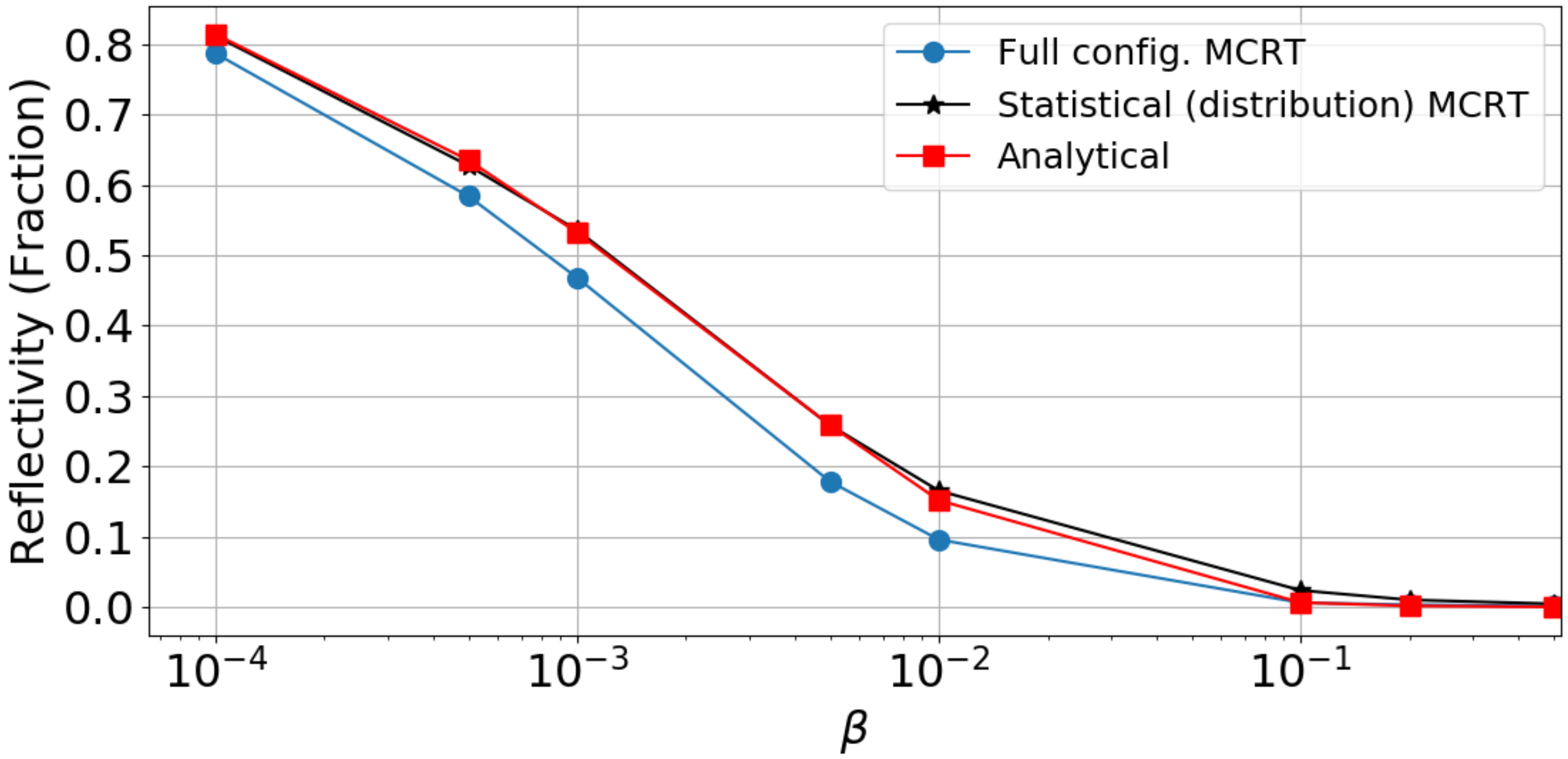

**Figure 11**. Reflectivity estimates from MCRT, pack-free MC, and purely analytical model using Theorem 2 with the numerical fitted distributions for the porous medium of Figure 8.

Finally, Fig. 12 compares the analytical reflectivity formula of Theorem 3 with exhaustive MC simulations for a non-Beerian two-sided porous medium with an exponential ($\mu$) mean-free beam length distribution for zero incident angle. The curve indicates a precise match, therefore validating the derivation.

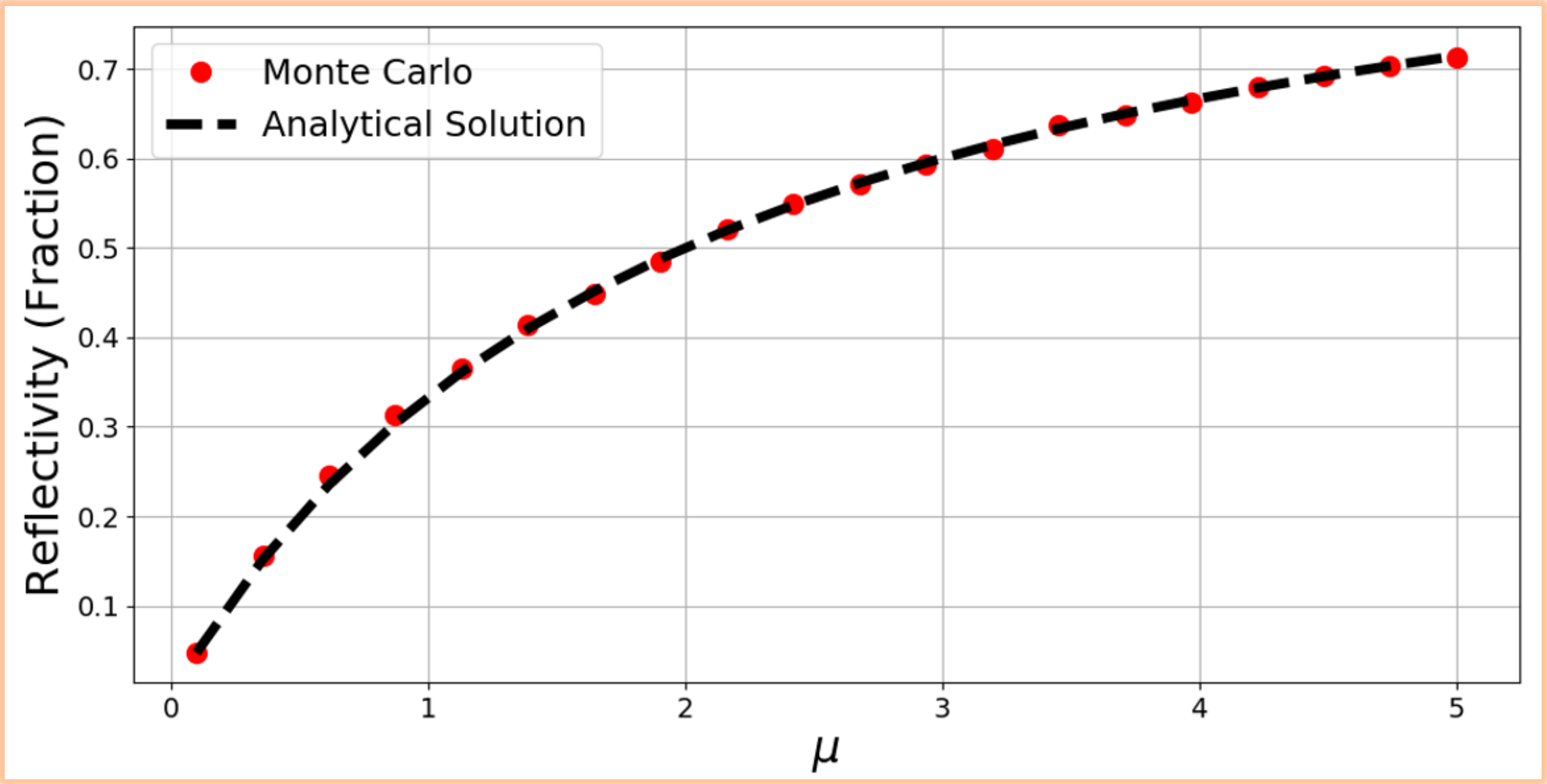

**Figure 12.** Comparison of the analytical reflectivity formula of Theorem 3 with exhaustive MC simulations for a non-Beerian two-sided porous medium with an exponential ($\mu$) mean-free beam length distribution.

## 4. Conclusion & Glimpse of Future Work

We extracted the existing machinery of Renewal, Ruin (Cramér-Lundberg), and surplus risk theories and used them to obtain preliminary precise geometric optics radiation estimations in porous media. Future work shall seek to expand Theorems 2&3 for probability decompositions that are more complex than a single rate exponential. Erlang mode decomposition can capture the complexity of arbitrary heterogeneous geometric factors. Future work shall strive to derive analytically simplified moment-generating function expressions and use techniques similar to those of this work to approximate spectral power reflection as a function of mean free beam length mode decomposition. As suggested by more recent renewal model literature[119], the overshoots and stopping times can be modeled by mixtures of Erlangs $\sum_{i=1}^{\infty} \pi'_i n_{k'_i,\mu'_i} x^{k'_i-1} e^{-x\mu'_i}$ when the step size distribution is also an Erlang mixture $\sum_{i=1}^{\infty} \pi_i n_{k_i,\mu_i} x^{k_i-1} e^{-x\mu_i}$. The literature machineries[114], [119] provide tools for the mapping of model mixture parameters $(\pi_i, k_i, \mu_i)$ to the solution mixture parameters $(\pi'_i, k'_i, \mu'_i)$ which, in the most complex case, involves finding roots of a polynomial equation of the order of the mixture size. Our future work shall strive to compile a thorough portfolio of derivations and possible closed forms, particularly in asymptotic limit regimes of parameter ratios using renewal processes[120]–[123]. Several other key results of Renewal theory can precisely explain asymptotic power behavior, uncertainty estimates of the derived calculations, and exact distribution functions. The work was limited to porous media with opaque particles. Future work shall consider refraction and other scattering modes. More advanced future work shall also consider two-dimensional renewal processes allowing for explicit modeling of other media shapes and more advanced parameters.

## References

[1] C. Argento and D. Bouvard, “A ray tracing method for evaluating the radiative heat transfer in porous media,” *International Journal of Heat and Mass Transfer*, vol. 39, no. 15, pp. 3175–3180, Oct. 1996, doi: 10.1016/0017-9310(95)00403-3.

[2] Y. S. Yang, J. R. Howell, and D. E. Klein, “Radiative Heat Transfer Through a Randomly Packed Bed of Spheres by the Monte Carlo Method,” *Journal of Heat Transfer*, vol. 105, no. 2, pp. 325–332, May 1983, doi: 10.1115/1.3245582.

[3] Y. Hua, G. Flamant, J. Lu, and D. Gauthier, “3D modelling of radiative heat transfer in circulating fluidized bed combustors: influence of the particulate composition,” *International Journal of Heat and Mass Transfer*, vol. 48, no. 6, pp. 1145–1154, Mar. 2005, doi: 10.1016/j.ijheatmasstransfer.2004.10.001.

[4] K. Kamiuto and S. San Yee, “Correlated radiative transfer through a packed bed of opaque spheres,” *International Communications in Heat and Mass Transfer*, vol. 32, no. 1–2, pp. 133–139, Jan. 2005, doi: 10.1016/j.icheatmasstransfer.2004.03.021.

[5] P. Rubiolo and J.-M. Gatt, “Modeling of the radiative contribution to heat transfer in porous media composed of spheres or cylinders,” *International Journal of Thermal Sciences*, vol. 41, no. 5, pp. 401–411, Apr. 2002, doi: 10.1016/S1290-0729(02)01332-7.

[6] P. von Zedtwitz, W. Lipiński, and A. Steinfeld, “Numerical and experimental study of gas–particle radiative heat exchange in a fluidized-bed reactor for steam-gasification of coal,” *Chemical Engineering Science*, vol. 62, no. 1–2, pp. 599–607, Jan. 2007, doi: 10.1016/j.ces.2006.09.027.

[7] J. R. Mahan, *The Monte Carlo ray-trace method in radiation heat transfer and applied optics*. Hoboken, NJ: John Wiley & Sons, 2018.

[8] P. M. Campbell, “Monte carlo method for radiative transfer,” *International Journal of Heat and Mass Transfer*, vol. 10, no. 4, pp. 519–527, Apr. 1967, doi: 10.1016/0017-9310(67)90171-8.

[9] J. H. Halton, “A Retrospective and Prospective Survey of the Monte Carlo Method,” *SIAM Rev.*, vol. 12, no. 1, pp. 1–63, Jan. 1970, doi: 10.1137/1012001.

[10] J. C. Chen and S. W. Churchill, “Radiant heat transfer in packed beds,” *AIChE J.*, vol. 9, no. 1, pp. 35–41, Jan. 1963, doi: 10.1002/aic.690090108.

[11] C.-L. Tien and K. Vafai, “Convective and radiative heat transfer in porous media,” *Advances in applied mechanics*, vol. 27, pp. 225–281, 1989.

[12] R. J. Cimini and J. C. Chen, “EXPERIMENTAL MEASUREMENTS OF RADIANT TRANSMISSION THROUGH PACKED AND FLUIDIZED MEDIA,” *Experimental Heat Transfer*, vol. 1, no. 1, pp. 45–56, Jan. 1987, doi: 10.1080/08916158708946330.

[13] D. Baillis and J.-F. Sacadura, “Thermal radiation properties of dispersed media: theoretical prediction and experimental characterization,” *Journal of Quantitative Spectroscopy and Radiative Transfer*, vol. 67, no. 5, pp. 327–363, Dec. 2000, doi: 10.1016/S0022-4073(99)00234-4.

[14] A. Menart, H. S. Lee, and R. O. Buckius, “EXPERIMENTAL DETERMINATION OF RADIATIVE PROPERTIES FOR SCATTERING PARTICULATE,” *Experimental Heat Transfer*, vol. 2, no. 4, pp. 309–332, Jan. 1989, doi: 10.1080/08916158908946371.

[15] A. É. Bugrov *et al.*, “Absorption and scattering of high-power laser radiation in low-density porous media,” *J. Exp. Theor. Phys.*, vol. 88, no. 3, pp. 441–448, Mar. 1999, doi: 10.1134/1.558814.

[16] Z. Sun and K. P. Shine, “Studies of the radiative properties of ice and mixed-phase clouds,” *Q.J Royal Met. Soc.*, vol. 120, no. 515, pp. 111–137, Jan. 1994, doi: 10.1002/qj.49712051508.
[17] W. B. Argo and J. M. Smith, “Heat transfer in packed beds-prediction of radial rates in gas-solid beds,” *Chemical Engineering Progress*, vol. 49, no. 8, pp. 443–451, 1953.
[18] C. K. Chan and C. L. Tien, “Radiative Transfer in Packed Spheres,” *Journal of Heat Transfer*, vol. 96, no. 1, pp. 52–58, Feb. 1974, doi: 10.1115/1.3450140.
[19] A. P. de Wasch and G. F. Froment, “Heat transfer in packed beds,” *Chemical Engineering Science*, vol. 27, no. 3, pp. 567–576, Mar. 1972, doi: 10.1016/0009-2509(72)87012-X.
[20] J. Randrianalisoa and D. Baillis, “Analytical model of radiative properties of packed beds and dispersed media,” *International Journal of Heat and Mass Transfer*, vol. 70, pp. 264–275, Mar. 2014, doi: 10.1016/j.ijheatmasstransfer.2013.10.071.
[21] L. Pilon* and R. Viskanta, “Radiation Characteristics of Glass Containing Gas Bubbles,” *Journal of the American Ceramic Society*, vol. 86, no. 8, pp. 1313–1320, Aug. 2003, doi: 10.1111/j.1151-2916.2003.tb03468.x.
[22] J. Randrianalisoa, D. Baillis, and L. Pilon, “Modeling radiation characteristics of semitransparent media containing bubbles or particles,” *J. Opt. Soc. Am. A*, vol. 23, no. 7, p. 1645, Jul. 2006, doi: 10.1364/JOSAA.23.001645.
[23] E. F. M. Van Der Held, “The contribution of radiation to the conduction of heat: II. Boundary conditions,” *Appl. sci. Res.*, vol. 4, no. 2, pp. 77–99, Mar. 1953, doi: 10.1007/BF03184939.
[24] Y. Zhao and G. H. Tang, “Monte Carlo study on extinction coefficient of silicon carbide porous media used for solar receiver,” *International Journal of Heat and Mass Transfer*, vol. 92, pp. 1061–1065, Jan. 2016, doi: 10.1016/j.ijheatmasstransfer.2015.08.105.
[25] C.-A. Wang, L.-X. Ma, J.-Y. Tan, and L.-H. Liu, “Study of radiative transfer in 1D densely packed bed layer containing absorbing–scattering spherical particles,” *International Journal of Heat and Mass Transfer*, vol. 102, pp. 669–678, Nov. 2016, doi: 10.1016/j.ijheatmasstransfer.2016.06.065.
[26] Y. Zhao, G. H. Tang, and M. Du, “Numerical study of radiative properties of nanoporous silica aerogel,” *International Journal of Thermal Sciences*, vol. 89, pp. 110–120, Mar. 2015, doi: 10.1016/j.ijthermalsci.2014.10.013.
[27] J. Randrianalisoa, S. Haussener, D. Baillis, and W. Lipiński, “Radiative characterization of random fibrous media with long cylindrical fibers: Comparison of single- and multi-RTE approaches,” *Journal of Quantitative Spectroscopy and Radiative Transfer*, vol. 202, pp. 220–232, Nov. 2017, doi: 10.1016/j.jqsrt.2017.08.002.
[28] B. X. Wang and C. Y. Zhao, “Modeling radiative properties of air plasma sprayed thermal barrier coatings in the dependent scattering regime,” *International Journal of Heat and Mass Transfer*, vol. 89, pp. 920–928, Oct. 2015, doi: 10.1016/j.ijheatmasstransfer.2015.06.017.
[29] S. Subramaniam and M. P. Mengüç, “Solution of the inverse radiation problem for inhomogeneous and anisotropically scattering media using a Monte Carlo technique,” *International Journal of Heat and Mass Transfer*, vol. 34, no. 1, pp. 253–266, Jan. 1991, doi: 10.1016/0017-9310(91)90192-H.
[30] S. Le Foll, F. André, A. Delmas, J. M. Bouilly, and Y. Aspa, “Radiative transfer modelling inside thermal protection system using hybrid homogenization method for a

backward Monte Carlo method coupled with Mie theory," *J. Phys.: Conf. Ser.*, vol. 369, p. 012028, Jun. 2012, doi: 10.1088/1742-6596/369/1/012028.
[31] M. Oguma, "Solution of two-dimensional blackbody inverse radiation problem by inverse Monte Carlo method," Maui, Hawaii, 1995.
[32] H. M. Park and M. C. Sung, "Solution of a multidimensional inverse radiation problem by means of mode reduction," *Int. J. Numer. Meth. Engng.*, vol. 60, no. 12, pp. 1949–1968, Jul. 2004, doi: 10.1002/nme.1028.
[33] S. W. Baek, D. Y. Byun, and S. J. Kang, "The combined Monte-Carlo and finite-volume method for radiation in a two-dimensional irregular geometry," *International Journal of Heat and Mass Transfer*, vol. 43, no. 13, pp. 2337–2344, Jul. 2000, doi: 10.1016/S0017-9310(99)00288-4.
[34] A. Safavinejad, S. H. Mansouri, and S. M. H. Sarvari, "Inverse boundary design of two-dimensional radiant enclosures with absorbing—emitting media using micro-genetic algorithm," *Proceedings of the Institution of Mechanical Engineers, Part C: Journal of Mechanical Engineering Science*, vol. 221, no. 8, pp. 945–948, Aug. 2007, doi: 10.1243/09544062JMES154.
[35] A. Safavinejad, S. H. Mansouri, A. Sakurai, and S. Maruyama, "Optimal number and location of heaters in 2-D radiant enclosures composed of specular and diffuse surfaces using micro-genetic algorithm," *Applied Thermal Engineering*, vol. 29, no. 5–6, pp. 1075–1085, Apr. 2009, doi: 10.1016/j.applthermaleng.2008.05.025.
[36] M. Mosavati, F. Kowsary, and B. Mosavati, "A Novel, Noniterative Inverse Boundary Design Regularized Solution Technique Using the Backward Monte Carlo Method," *Journal of Heat Transfer*, vol. 135, no. 4, p. 042701, Apr. 2013, doi: 10.1115/1.4022994.
[37] B. Mosavati, M. Mosavati, and F. Kowsary, "Solution of radiative inverse boundary design problem in a combined radiating-free convecting furnace," *International Communications in Heat and Mass Transfer*, vol. 45, pp. 130–136, Jul. 2013, doi: 10.1016/j.icheatmasstransfer.2013.04.011.
[38] R. B. Mulford, N. S. Collins, M. S. Farnsworth, M. R. Jones, and B. D. Iverson, "Total hemispherical apparent radiative properties of the infinite V-groove with specular reflection," *International Journal of Heat and Mass Transfer*, vol. 124, pp. 168–176, Sep. 2018, doi: 10.1016/j.ijheatmasstransfer.2018.03.041.
[39] M. Yarahmadi, J. R. Mahan, and F. Kowsary, "A New Approach to Inverse Boundary Design in Radiation Heat Transfer," in *Advances in Heat Transfer and Thermal Engineering*, 2021, pp. 377–383.
[40] L. Ibarrart *et al.*, "COMBINED CONDUCTIVE-CONVECTIVE-RADIATIVE HEAT TRANSFER IN COMPLEX GEOMETRY USING THE MONTE CARLO METHOD : APPLICATION TO SOLAR RECEIVERS," in *International Heat Transfer Conference 16*, Beijing, China, 2018, pp. 8135–8142. doi: 10.1615/IHTC16.pma.023662.
[41] D. G. Coronell and K. F. Jensen, "A Monte Carlo Simulation Study of Radiation Heat Transfer in the Multiwafer LPCVD Reactor," *J. Electrochem. Soc.*, vol. 141, no. 2, pp. 496–501, Feb. 1994, doi: 10.1149/1.2054753.
[42] F. Wang, Y. Shuai, H. Tan, and C. Yu, "Thermal performance analysis of porous media receiver with concentrated solar irradiation," *International Journal of Heat and Mass Transfer*, vol. 62, pp. 247–254, Jul. 2013, doi: 10.1016/j.ijheatmasstransfer.2013.03.003.

[43] Y. Zhao and G. H. Tang, "Monte Carlo Study on Carbon-Gradient-Doped Silica Aerogel Insulation," *j nanosci nanotechnol*, vol. 15, no. 4, pp. 3259–3264, Apr. 2015, doi: 10.1166/jnn.2015.9669.

[44] X. Ji, Y. Zhang, X. Li, G. Fan, and M. Li, "Solar ray collection rate fluctuation analysis with Monte Carlo Ray Tracing method for space solar power satellite," *Solar Energy*, vol. 185, pp. 235–244, Jun. 2019, doi: 10.1016/j.solener.2019.04.067.

[45] M. Rafiee, S. Chandra, H. Ahmed, K. Barnham, and S. J. McCormack, "Optical Coupling Sensitivity Study of Luminescent PV Devices Using Monte Carlo Ray-Tracing Model," in *Renewable Energy and Sustainable Buildings*, Cham: Springer, 2020, pp. 869–877.

[46] M. Rafiee, H. Ahmed, S. Chandra, A. Sethi, and S. J. McCormack, "Monte Carlo Ray Tracing Modelling of Multi-Crystalline Silicon Photovoltaic Device Enhanced by Luminescent Material," in *2018 IEEE 7th World Conference on Photovoltaic Energy Conversion (WCPEC) (A Joint Conference of 45th IEEE PVSC, 28th PVSEC & 34th EU PVSEC)*, Waikoloa Village, HI, Jun. 2018, pp. 3139–3141. doi: 10.1109/PVSC.2018.8547436.

[47] P. Parthasarathy, P. Habisreuther, and N. Zarzalis, "Identification of radiative properties of reticulated ceramic porous inert media using ray tracing technique," *Journal of Quantitative Spectroscopy and Radiative Transfer*, vol. 113, no. 15, pp. 1961–1969, Oct. 2012, doi: 10.1016/j.jqsrt.2012.05.017.

[48] D. Mackay and J. Cameron, "Introduction to monte carlo methods," in *Learning in graphical models*, Dordrecht: Springer, 1998, pp. 175–204.

[49] W. R. Gilks, S. Richardson, and D. Spiegelhalter, Eds., *Markov Chain Monte Carlo in Practice*, 0 ed. Chapman and Hall/CRC, 1995. doi: 10.1201/b14835.

[50] C. P. Robert and G. Casella, *Monte Carlo statistical methods*, 2nd ed. New York: Springer, 2004.

[51] M. Q. Brewster and C. L. Tien, "Radiative Transfer in Packed Fluidized Beds: Dependent Versus Independent Scattering," *Journal of Heat Transfer*, vol. 104, no. 4, pp. 573–579, Nov. 1982, doi: 10.1115/1.3245170.

[52] J. D. Cartigny, Y. Yamada, and C. L. Tien, "Radiative Transfer With Dependent Scattering by Particles: Part 1—Theoretical Investigation," *Journal of Heat Transfer*, vol. 108, no. 3, pp. 608–613, Aug. 1986, doi: 10.1115/1.3246979.

[53] Y. Yamada, J. D. Cartigny, and C. L. Tien, "Radiative Transfer With Dependent Scattering by Particles: Part 2—Experimental Investigation," *Journal of Heat Transfer*, vol. 108, no. 3, pp. 614–618, Aug. 1986, doi: 10.1115/1.3246980.

[54] C. L. Tien, "Thermal Radiation in Packed and Fluidized Beds," *Journal of Heat Transfer*, vol. 110, no. 4b, pp. 1230–1242, Nov. 1988, doi: 10.1115/1.3250623.

[55] B. L. Drolen and C. L. Tien, "Independent and dependent scattering in packed-sphere systems," *Journal of Thermophysics and Heat Transfer*, vol. 1, no. 1, pp. 63–68, Jan. 1987, doi: 10.2514/3.8.

[56] S. Kumar and C. L. Tien, "Dependent Absorption and Extinction of Radiation by Small Particles," *Journal of Heat Transfer*, vol. 112, no. 1, pp. 178–185, Feb. 1990, doi: 10.1115/1.2910342.

[57] K. Kamiuto, "Correlated radiative transfer in packed-sphere systems," *Journal of Quantitative Spectroscopy and Radiative Transfer*, vol. 43, no. 1, pp. 39–43, Jan. 1990, doi: 10.1016/0022-4073(90)90063-C.

[58] B. P. Singh and M. Kaviany, “Independent theory versus direct simulation of radiation heat transfer in packed beds,” *International Journal of Heat and Mass Transfer*, vol. 34, no. 11, pp. 2869–2882, Nov. 1991, doi: 10.1016/0017-9310(91)90247-C.
[59] B. P. Singh and M. Kaviany, “Modelling radiative heat transfer in packed beds,” *International Journal of Heat and Mass Transfer*, vol. 35, no. 6, pp. 1397–1405, Jun. 1992, doi: 10.1016/0017-9310(92)90031-M.
[60] R. Siegel and J. R. Howell, *Thermal radiation heat transfer*, 4th ed. New York: Taylor & Francis, 2002.
[61] L. A. Dombrovskij and L. A. Dombrovsky, *Radiation heat transfer in disperse systems*. New York: Begell House, 1996.
[62] E. M. Sparrow, *Radiation Heat Transfer, Augmented Edition*. 2017. Accessed: Jul. 17, 2022. [Online]. Available: https://www.vlebooks.com/vleweb/product/openreader?id=none&isbn=9781351420105
[63] Y. T. Feng and K. Han, “An accurate evaluation of geometric view factors for modelling radiative heat transfer in randomly packed beds of equally sized spheres,” *International Journal of Heat and Mass Transfer*, vol. 55, no. 23–24, pp. 6374–6383, Nov. 2012, doi: 10.1016/j.ijheatmasstransfer.2012.06.025.
[64] S. T. Flock, M. S. Patterson, B. C. Wilson, and D. R. Wyman, “Monte Carlo modeling of light propagation in highly scattering tissues. I. Model predictions and comparison with diffusion theory,” *IEEE Trans. Biomed. Eng.*, vol. 36, no. 12, pp. 1162–1168, Dec. 1989, doi: 10.1109/TBME.1989.1173624.
[65] É. Sanchez-Palencia, *Non-homogeneous media and vibration theory*. Berlin: Springer, 1980.
[66] E. Pasternak and H.-B. Mühlhaus, “Generalised homogenisation procedures for granular materials,” *J Eng Math*, vol. 52, no. 1, pp. 199–229, Jul. 2005, doi: 10.1007/s10665-004-3950-z.
[67] G. A. Pavliotis and A. M. Stuart, *Multiscale methods: averaging and homogenization*. New York: Springer, 2008.
[68] C. M. Rooney, C. P. Please, and S. D. Howison, “Homogenisation applied to thermal radiation in porous media,” *Eur. J. Appl. Math*, vol. 32, no. 5, pp. 784–805, Oct. 2021, doi: 10.1017/S0956792520000388.
[69] Q. Zhou, H. W. Zhang, and Y. G. Zheng, “A homogenization technique for heat transfer in periodic granular materials,” *Advanced Powder Technology*, vol. 23, no. 1, pp. 104–114, Jan. 2012, doi: 10.1016/j.apt.2011.01.002.
[70] Q. Brewster, “Volume Scattering of Radiation in Packed Beds of Large, Opaque Spheres,” *Journal of Heat Transfer*, vol. 126, no. 6, pp. 1048–1050, Dec. 2004, doi: 10.1115/1.1795247.
[71] M. F. Modest and S. Mazumder, *Radiative heat transfer*. Academic press, 2021.
[72] H. C. Hottel and A. F. Sarofim, *Radiative Transfer*. New York: McGraw-Hill, 1967.
[73] J. R. Howell, *A catalog of radiation configuration factors*. New York: McGraw-Hill Book Co, 1982.
[74] B.-H. Gao, H. Qi, D.-H. Jiang, Y.-T. Ren, and M.-J. He, “Efficient equation-solving integral equation method based on the radiation distribution factor for calculating radiative transfer in 3D anisotropic scattering medium,” *Journal of Quantitative Spectroscopy and Radiative Transfer*, vol. 275, p. 107886, Nov. 2021, doi: 10.1016/j.jqsrt.2021.107886.

[75] de C. Augusto, L. Dobrowolski, B. Giacomet, and N. Mendes, "Numerical method for calculating view factor between two surfaces," 2007.
[76] R. Coquard and D. Baillis, "Radiative Characteristics of Opaque Spherical Particles Beds: A New Method of Prediction," *Journal of Thermophysics and Heat Transfer*, vol. 18, no. 2, pp. 178–186, Apr. 2004, doi: 10.2514/1.5082.
[77] T. J. Hendricks and J. R. Howell, "Absorption/Scattering Coefficients and Scattering Phase Functions in Reticulated Porous Ceramics," *Journal of Heat Transfer*, vol. 118, no. 1, pp. 79–87, Feb. 1996, doi: 10.1115/1.2824071.
[78] M. Tancrez and J. Taine, "Direct identification of absorption and scattering coefficients and phase function of a porous medium by a Monte Carlo technique," *International Journal of Heat and Mass Transfer*, vol. 47, no. 2, pp. 373–383, Jan. 2004, doi: 10.1016/S0017-9310(03)00146-7.
[79] J. Petrasch, P. Wyss, and A. Steinfeld, "Tomography-based Monte Carlo determination of radiative properties of reticulate porous ceramics," *Journal of Quantitative Spectroscopy and Radiative Transfer*, vol. 105, no. 2, pp. 180–197, Jun. 2007, doi: 10.1016/j.jqsrt.2006.11.002.
[80] S. Haussener, P. Coray, W. Lipiński, P. Wyss, and A. Steinfeld, "Tomography-Based Heat and Mass Transfer Characterization of Reticulate Porous Ceramics for High-Temperature Processing," *Journal of Heat Transfer*, vol. 132, no. 2, p. 023305, Feb. 2010, doi: 10.1115/1.4000226.
[81] B. Zeghondy, E. Iacona, and J. Taine, "Experimental and RDFI calculated radiative properties of a mullite foam," *International Journal of Heat and Mass Transfer*, vol. 49, no. 19–20, pp. 3702–3707, Sep. 2006, doi: 10.1016/j.ijheatmasstransfer.2006.02.036.
[82] M. Zarrouati, F. Enguehard, and J. Taine, "Statistical characterization of near-wall radiative properties of a statistically non-homogeneous and anisotropic porous medium," *International Journal of Heat and Mass Transfer*, vol. 67, pp. 776–783, Dec. 2013, doi: 10.1016/j.ijheatmasstransfer.2013.08.021.
[83] B. X. Wang and C. Y. Zhao, "Effect of anisotropy on thermal radiation transport in porous ceramics," *International Journal of Thermal Sciences*, vol. 111, pp. 301–309, Jan. 2017, doi: 10.1016/j.ijthermalsci.2016.09.012.
[84] J. Taine and E. Iacona, "Upscaling Statistical Methodology for Radiative Transfer in Porous Media: New Trends," *Journal of Heat Transfer*, vol. 134, no. 3, p. 031012, Mar. 2012, doi: 10.1115/1.4005133.
[85] J. Taine, E. Iacona, and F. Bellet, "Radiation in porous media: an upscaling methodology," 2008, p. ISBN-978.
[86] S. García-Pareja, A. M. Lallena, and F. Salvat, "Variance-Reduction Methods for Monte Carlo Simulation of Radiation Transport," *Front. Phys.*, vol. 9, p. 718873, Oct. 2021, doi: 10.3389/fphy.2021.718873.
[87] G. Rubino and B. Tuffin, *Rare event simulation using Monte Carlo methods*. Chichester, U.K: Wiley, 2009.
[88] E. Woodcock, T. Murphy, P. Hemmings, and S. Longworth, "Techniques used in the GEM code for Monte Carlo neutronics calculations in reactors and other systems of complex geometry," in *Proc. Conf. Applications of Computing Methods to Reactor Problems*, 1965, vol. 557 no. 2.
[89] S. M. Kersch Alfred, "A FAST MONTE CARLO SCHEME FOR THERMAL RADIATION IN SEMICONDUCTOR PROCESSING APPLICATIONS," *Numerical*

*Heat Transfer, Part B: Fundamentals*, vol. 37, no. 2, pp. 185–199, Mar. 2000, doi: 10.1080/104077900275486.
[90] "MULTIPLE-RAYS TRACING TECHNIQUE FOR RADIATIVE EXCHANGE WITHIN PACKED BEDS," *Numerical Heat Transfer, Part B: Fundamentals*, vol. 37, no. 4, pp. 469–487, Jun. 2000, doi: 10.1080/10407790050051155.
[91] W. A. Coleman, "Mathematical Verification of a Certain Monte Carlo Sampling Technique and Applications of the Technique to Radiation Transport Problems," *Nuclear Science and Engineering*, vol. 32, no. 1, pp. 76–81, Apr. 1968, doi: 10.13182/NSE68-1.
[92] W. Lipiński, J. Petrasch, and S. Haussener, "Application of the spatial averaging theorem to radiative heat transfer in two-phase media," *Journal of Quantitative Spectroscopy and Radiative Transfer*, vol. 111, no. 1, pp. 253–258, Jan. 2010, doi: 10.1016/j.jqsrt.2009.08.001.
[93] J. Petrasch, S. Haussener, and W. Lipiński, "Discrete vs. continuum-scale simulation of radiative transfer in semitransparent two-phase media," *Journal of Quantitative Spectroscopy and Radiative Transfer*, vol. 112, no. 9, pp. 1450–1459, Jun. 2011, doi: 10.1016/j.jqsrt.2011.01.025.
[94] A. V. Gusarov, "Homogenization of radiation transfer in two-phase media with irregular phase boundaries," *Phys. Rev. B*, vol. 77, no. 14, p. 144201, Apr. 2008, doi: 10.1103/PhysRevB.77.144201.
[95] R. Coquard and D. Baillis, "Radiative characteristics of beds made of large spheres containing an absorbing and scattering medium," *International Journal of Thermal Sciences*, vol. 44, no. 10, pp. 926–932, Oct. 2005, doi: 10.1016/j.ijthermalsci.2005.03.009.
[96] R. Coquard and D. Baillis, "Radiative Characteristics of Beds of Spheres Containing an Absorbing and Scattering Medium.," *Journal of Thermophysics and Heat Transfer*, vol. 19, no. 2, pp. 226–234, Apr. 2005, doi: 10.2514/1.6809.
[97] D. Moser, S. Pannala, and J. Murthy, "Computation of Effective Radiative Properties of Powders for Selective Laser Sintering Simulations," *JOM*, vol. 67, no. 5, pp. 1194–1202, May 2015, doi: 10.1007/s11837-015-1386-8.
[98] M. Yarahmadi, J. Robert Mahan, and K. McFall, "Artificial Neural Networks in Radiation Heat Transfer Analysis," *Journal of Heat Transfer*, vol. 142, no. 9, p. 092801, Sep. 2020, doi: 10.1115/1.4047052.
[99] M. Taki, A. Rohani, and H. Yildizhan, "Application of machine learning for solar radiation modeling," *Theor Appl Climatol*, vol. 143, no. 3–4, pp. 1599–1613, Feb. 2021, doi: 10.1007/s00704-020-03484-x.
[100] S. Hajimirza and H. Sharadga, "Learning thermal radiative properties of porous media from engineered geometric features," *International Journal of Heat and Mass Transfer*, vol. 179, p. 121668, Nov. 2021, doi: 10.1016/j.ijheatmasstransfer.2021.121668.
[101] J. J. García-Esteban, J. Bravo-Abad, and J. C. Cuevas, "Deep Learning for the Modeling and Inverse Design of Radiative Heat Transfer," *Phys. Rev. Applied*, vol. 16, no. 6, p. 064006, Dec. 2021, doi: 10.1103/PhysRevApplied.16.064006.
[102] J. Tausendschön, G. Stöckl, and S. Radl, "Machine Learning for heat radiation modeling of bi- and polydisperse particle systems including walls," *Particuology*, vol. 74, pp. 119–140, Mar. 2023, doi: 10.1016/j.partic.2022.05.011.
[103] Q. Zou, N. Gui, X. Yang, J. Tu, and S. Jiang, "Comparative Study on the Numerical Methods for View Factor Computation for Packed Pebble Beds: Back Propagation Neural

Network Methods Versus Monte Carlo Methods," *Journal of Heat Transfer*, vol. 143, no. 8, p. 083301, Aug. 2021, doi: 10.1115/1.4051075.
[104] W. W. Yuen, C. L. Chow, and W. C. Tam, "Analysis of Radiative Heat Transfer in Inhomogeneous Nonisothermal Media Using Neural Networks," *Journal of Thermophysics and Heat Transfer*, vol. 30, no. 4, pp. 897–911, Oct. 2016, doi: 10.2514/1.T4805.
[105] M. Chevrollier, "Radiation trapping and Lévy flights in atomic vapours: an introductory review," *Contemporary Physics*, vol. 53, no. 3, pp. 227–239, May 2012, doi: 10.1080/00107514.2012.684481.
[106] L. Fredriksson, "A brief survey of Lévy walks: with applications to probe diffusion," 2010.
[107] M. F. Shlesinger, J. Klafter, and B. J. West, "Levy walks with applications to turbulence and chaos," *Physica A: Statistical Mechanics and its Applications*, vol. 140, no. 1–2, pp. 212–218, Dec. 1986, doi: 10.1016/0378-4371(86)90224-4.
[108] B. Berkowitz, A. Cortis, M. Dentz, and H. Scher, "Modeling non-Fickian transport in geological formations as a continuous time random walk," *Rev. Geophys.*, vol. 44, no. 2, p. RG2003, 2006, doi: 10.1029/2005RG000178.
[109] I. Alemany, J. N. Rose, J. Garnier-Brun, A. D. Scott, and D. J. Doorly, "Random walk diffusion simulations in semi-permeable layered media with varying diffusivity," *Sci Rep*, vol. 12, no. 1, p. 10759, Dec. 2022, doi: 10.1038/s41598-022-14541-y.
[110] T.-C. Yeh, R. Khaleel, and K. C. Carroll, *Flow through heterogeneous geologic media*. Cambridge University Press, 2015.
[111] A. Wagner *et al.*, "Permeability Estimation of Regular Porous Structures: A Benchmark for Comparison of Methods," *Transp Porous Med*, vol. 138, no. 1, pp. 1–23, May 2021, doi: 10.1007/s11242-021-01586-2.
[112] C.-O. Hwang, J. A. Given, and M. Mascagni, "On the rapid estimation of permeability for porous media using Brownian motion paths," *Physics of Fluids*, vol. 12, no. 7, pp. 1699–1709, Jul. 2000, doi: 10.1063/1.870420.
[113] N. A. Simonov and M. Mascagni, "Random Walk Algorithms for Estimating Effective Properties of Digitized Porous Media *," *Monte Carlo Methods and Applications*, vol. 10, no. 3–4, Jan. 2004, doi: 10.1515/mcma.2004.10.3-4.599.
[114] M. Jacobsen, "Exit times for a class of random walks exact distribution results," *Journal of Applied Probability*, vol. 48, no. A, pp. 51–63, Aug. 2011, doi: 10.1239/jap/1318940455.
[115] M. Dekking, Ed., *A modern introduction to probability and statistics: understanding why and how*. London: Springer, 2005.
[116] P. S. Bullen, *Dictionary of inequalities*, 2. ed. Boca Raton, Fla.: CRC Press, 2015.
[117] C. Cox, *Delta method*, vol. 2. 2005.
[118] J. Janssen and R. Manca, *Applied semi-Markov processes*. New York, NY: Springer Science+Business Media, 2006.
[119] M. V. Boutsikas and K. Politis, "Exit Times, Overshoot and Undershoot for a Surplus Process in the Presence of an Upper Barrier," *Methodol Comput Appl Probab*, vol. 19, no. 1, pp. 75–95, Mar. 2017, doi: 10.1007/s11009-015-9459-2.
[120] D. C. M. Dickson, "On the distribution of the surplus prior to ruin," *Insurance: Mathematics and Economics*, vol. 11, no. 3, pp. 191–207, Oct. 1992, doi: 10.1016/0167-6687(92)90026-8.

[121] D. C. M. Dickson, “On a Class of Renewal Risk Processes,” *North American Actuarial Journal*, vol. 2, no. 3, pp. 60–68, Jul. 1998, doi: 10.1080/10920277.1998.10595723.
[122] D. C. M. Dickson and J. R. Gray, “Approximations to ruin probability in the presence of an upper absorbing barrier,” *Scandinavian Actuarial Journal*, vol. 1984, no. 2, pp. 105–115, Apr. 1984, doi: 10.1080/03461238.1984.10413758.
[123] J.-M. Reinhard and M. Snoussi, “On the distribution of the surplus prior to ruin in a discrete semi-Markov risk model,” *ASTIN Bulletin: The Journal of the IAA 31*, vol. 2, pp. 255–273, 2001.
[124] A. Gut, *Stopped Random Walks*. New York, NY: Springer New York, 2009. doi: 10.1007/978-0-387-87835-5.
[125] T. L. Lai, “Asymptotic Moments of Random Walks with Applications to Ladder Variables and Renewal Theory,” *Ann. Probab.*, vol. 4, no. 1, Feb. 1976, doi: 10.1214/aop/1176996180.
[126] G. Englund, “Remainder term estimate for the asymptotic normality of the number of renewals,” *Journal of Applied Probability*, vol. 17, no. 4, pp. 1108–1113, Dec. 1980, doi: 10.2307/3213222.
[127] J. L. Teugels, “Renewal Theorems When the First or the Second Moment is Infinite,” *Ann. Math. Statist.*, vol. 39, no. 4, pp. 1210–1219, Aug. 1968, doi: 10.1214/aoms/1177698246.
[128] N. R. Mohan, “Teugels’ Renewal Theorem and Stable Laws,” *Ann. Probab.*, vol. 4, no. 5, Oct. 1976, doi: 10.1214/aop/1176995991.